\documentclass[12pt]{amsart}
\usepackage[margin=1.05in]{geometry}
\usepackage{amsmath,amssymb,amsthm,mathtools}
\usepackage{xcolor}
\usepackage{booktabs,array}
\usepackage{tikz}
\usepackage{float}
\usetikzlibrary{decorations.pathreplacing,arrows.meta,patterns}
\usepackage{pgfplots}
\pgfplotsset{compat=1.16}
\definecolor{lava}{rgb}{0.81,0.06,0.13}
\definecolor{Cblue}{rgb}{0.50,0.85,0.85}
\definecolor{lime}{HTML}{A6CE39}
\definecolor{figA}{HTML}{2A78D6}
\definecolor{figB}{HTML}{EB6834}
\definecolor{figC}{HTML}{1BAF7A}
\definecolor{figGrid}{HTML}{BBBBBB}

\usepackage[colorlinks,linkcolor=lava,citecolor=blue,urlcolor=brown,hypertexnames=false]{hyperref}

\newtheorem{theorem}{Theorem}[section]
\newtheorem{proposition}[theorem]{Proposition}
\newtheorem{lemma}[theorem]{Lemma}
\newtheorem{corollary}[theorem]{Corollary}
\newtheorem{conjecture}[theorem]{Conjecture}
\theoremstyle{definition}
\newtheorem{definition}[theorem]{Definition}
\newtheorem{question}[theorem]{Question}
\theoremstyle{remark}
\newtheorem{remark}[theorem]{Remark}

\numberwithin{equation}{section}

\newcommand{\R}{\mathbb{R}}
\newcommand{\C}{\mathbb{C}}
\newcommand{\Sch}{\mathcal{S}(\R)}
\newcommand{\LS}{\mathrm{LS}}
\newcommand{\sech}{\operatorname{sech}}
\newcommand{\Hdot}{\dot H^{1/2}}
\newcommand{\eps}{\varepsilon}
\newcommand{\Mass}{M}
\newcommand{\Cm}{\mathcal{C}}
\newcommand{\sol}{\mathfrak{S}}

\title[Uniform local smoothing for DNLS]{Uniform local smoothing for the derivative NLS and the soliton mass threshold}

\author{Mohamed Majdoub}
\address{Department of Mathematics, College of Science, Imam Abdulrahman Bin Faisal University, P.O. Box 1982, Dammam, Saudi Arabia}
\address{Basic and Applied Scientific Research Center, Imam Abdulrahman Bin Faisal University, P.O. Box 1982, 31441, Dammam, Saudi Arabia}
\email{mmajdoub@iau.edu.sa}
\email{med.majdoub@gmail.com}
\email{mohamed.majdoub@fst.rnu.tn}
\urladdr{https://orcid.org/0000-0001-6038-1069}

\date{\today}

\begin{document}

\begin{abstract}
Harrop-Griffiths, Killip, Ntekoume and Vi\c{s}an proved that the derivative nonlinear Schr\"odinger equation (DNLS) is globally well-posed in $L^2(\R)$. A key ingredient of their proof is a local smoothing estimate. It holds on $L^2$-equicontinuous sets of initial data, but it cannot hold on $L^2$-balls, because rescaled stationary solitons, all of mass $2\pi$, violate it. This raises a natural question: for which masses $m$ does local smoothing hold uniformly on the ball $\{\Mass(q_0)\le m\}$? We prove three results. First, we give a quantitative version of the local smoothing estimate in which the data enter only through the mass and an explicit frequency scale $\kappa(q_0)$. On time intervals of length $2T$ the bound grows like $1+\kappa(q_0)^2T$, and the stationary soliton shows that the power of $\kappa$ cannot be lowered below one. Second, we prove uniform local smoothing on balls of small mass, together with a scale-invariant Morrey-type form and a localized $L^4_{t,x}$ bound. Third, we analyse the soliton family from both sides. The uniform constant must blow up at least like $(2\pi-m)^{-1}$ as $m\uparrow 2\pi$, while every soliton of mass $m<2\pi$ satisfies local smoothing with a constant $O\big((2\pi-m)^{-1}\big)$. These facts single out the stationary-soliton mass $2\pi$, which is also the Hayashi--Ozawa coercivity threshold, as the natural candidate for the sharp threshold. We formulate this as a conjecture, give a spectral reformulation in terms of the Kaup--Newell scattering data, and discuss the obstacles to a proof.
\end{abstract}

\maketitle

\setcounter{tocdepth}{1}

\section{Introduction and main results}\label{sec:intro}

We consider the derivative nonlinear Schr\"odinger equation on the line,
\begin{equation}\label{DNLS}
i\partial_t q + \partial_x^2 q + i\partial_x\big(|q|^2 q\big)=0,\qquad q:\R_t\times\R_x\to\C .
\end{equation}
It was derived as a model for circularly polarized Alfv\'en waves propagating along a magnetic field in a plasma \cite{MioEtAl76,Mjolhus76}, and it is completely integrable through the Kaup--Newell Lax pair \cite{KaupNewell}. The equation conserves the mass $\Mass(q)=\int_\R|q|^2\,dx$ and is invariant under the scaling
\begin{equation}\label{scaling}
q(t,x)\longmapsto q_\lambda(t,x):=\lambda^{1/2}q(\lambda^2t,\lambda x),\qquad \lambda>0,
\end{equation}
which leaves $\Mass$ unchanged. Hence $L^2(\R)$ is the scaling-critical space.

In a recent landmark paper, Harrop-Griffiths, Killip, Ntekoume and Vi\c{s}an \cite{HKNV} proved that \eqref{DNLS} is globally well-posed in $L^2(\R)$. More precisely, the data-to-solution map, defined on Schwartz data, extends to a jointly continuous map
\begin{equation}\label{Phi}
\Phi:\R\times L^2(\R)\to L^2(\R),\qquad q(t)=\Phi(t)q_0 .
\end{equation}
This result is the end point of a long line of work, which we recall in Section~\ref{subsec:lit}. Its proof relies on the second-generation method of commuting flows. The analytic heart of the method consists of two families of estimates. The first are local smoothing estimates, obtained from microscopic conservation laws for the logarithm of the Kaup--Newell perturbation determinant. The second express the tightness of orbits.

Let us recall the local smoothing statement of \cite{HKNV}. We use the localizing weights
\begin{equation}\label{weight}
\psi(x):=\big(\sech(x/99)\big)^{1/2},\qquad W:=\psi^{12}=\sech^6(\cdot/99),\qquad W_h:=W(\cdot-h).
\end{equation}
By \cite[Theorem~1.5]{HKNV}, if $Q\subset\Sch$ is bounded and equicontinuous in $L^2$, then for every $T>0$
\begin{equation}\label{HKNV-LS}
\sup_{h\in\R}\int_{-T}^{T}\big\|W_h\,q(t)\big\|_{H^{1/2}}^2\,dt\lesssim_{T,Q}\|q(0)\|_{L^2}^2,\qquad q(0)\in Q .
\end{equation}
The authors stress that the estimate is claimed for equicontinuous sets and not for balls, and that this restriction is necessary. Indeed, \cite[Proposition~1.11]{HKNV} exhibits solutions $q_n$, obtained by rescaling the stationary soliton, such that $\Mass(q_n)\equiv 2\pi$ while a localized local smoothing norm of $q_n$ on $[-1,1]$ tends to infinity. In \eqref{HKNV-LS}, the implicit constant comes from a rescaling that makes the relevant set of Schwartz functions ``$\delta$-good'' in the sense of \cite[Definition~1.8]{HKNV}. As the authors note in the proof of \cite[Theorem~1.5]{HKNV}, this constant depends on $Q$ only through the rescaling parameter, but that dependence is not made explicit.

The counterexample of \cite{HKNV} lives exactly on the sphere $\Mass=2\pi$, and $2\pi$ is a distinguished number for DNLS. It is the mass of the stationary soliton \cite[(1.25), (1.28)]{HKNV}. It is also the threshold below which the energy is coercive, a fact used by Hayashi and Ozawa \cite{HayashiOzawa92,HayashiOzawa94} to obtain global $H^1$ solutions. Moreover, Bahouri and Perelman point out \cite{BahouriPerelman22} that $2\pi$ and $4\pi$ are the masses of the static bright soliton and of the algebraic soliton, whose profiles are extremals of the relevant Gagliardo--Nirenberg inequalities. This suggests the following problem.

\medskip
\noindent\textbf{Problem.} \emph{Determine the largest $m^*\in(0,\infty]$ such that local smoothing holds uniformly on every ball $\{q_0\in L^2:\Mass(q_0)\le m\}$ with $m<m^*$, with a constant depending only on $m$. Describe the size of the optimal constant, and its dependence on the data when equicontinuity replaces smallness.}
\medskip

To make this precise, we introduce, for a solution $q$ of \eqref{DNLS} and $T>0$, the local smoothing functional
\begin{equation}\label{LS}
\LS_T(q):=\sup_{h\in\R}\int_{-T}^{T}\big\|W_h\,q(t)\big\|_{H^{1/2}}^2\,dt ,
\end{equation}
where $\|f\|_{H^{s}}^2=\int_\R(4+\xi^2)^{s}|\hat f(\xi)|^2\,d\xi$ as in \cite{HKNV}. We also introduce the \emph{uniform local smoothing constant}
\begin{equation}\label{Cm}
\Cm(m):=\sup\Big\{\frac{\LS_1\big(\Phi(\cdot)q_0\big)}{\Mass(q_0)}:\ q_0\in L^2(\R),\ 0<\Mass(q_0)\le m\Big\}\in(0,\infty] .
\end{equation}
For $q_0\in L^2$ the integrand in $\LS_1(\Phi(\cdot)q_0)$ is a lower semicontinuous, hence measurable, function of $t$, because $t\mapsto\Phi(t)q_0$ is continuous in $L^2$ and the $H^{1/2}$ norm is lower semicontinuous with respect to $L^2$ convergence (see the proof of Corollary~\ref{cor:small}). The function $\Cm$ is nondecreasing, and we set $m^*:=\sup\{m>0:\Cm(m)<\infty\}$. By \cite[Proposition~1.11]{HKNV}, or by Theorem~\ref{thm:B}(i) below, $m^*\le 2\pi$.

\subsection{Main results}

Throughout the paper we fix $0<\eps<1/2$ (for definiteness $\eps=1/4$; this parameter is denoted by $\sigma$ in \cite{HKNV}). We let $\delta_0>0$ be the smallness parameter for which \cite[Proposition~4.1]{HKNV} holds; see Proposition~\ref{prop:HKNV41} below. For $\kappa>0$ we use the equicontinuity quantity of \cite{HKNV},
\begin{equation}\label{Enorm}
\|f\|_{E^\eps_\kappa}^2:=\int_\R\frac{|\xi|^{2\eps}}{(4\kappa^2+\xi^2)^{\eps}}\,|\hat f(\xi)|^2\,d\xi ,
\end{equation}
which is the norm $\|f\|_{E^\sigma_{\sigma,\kappa}}$ of \cite[(2.2)]{HKNV} with $\sigma=\eps$; for $\kappa=1$ its square is the left-hand side of \cite[(1.19)]{HKNV}.

\begin{definition}[Frequency scale of an orbit]\label{def:kappa}
For $q_0\in\Sch$ with solution $q(t)=\Phi(t)q_0$ we set
\[
\kappa(q_0):=\inf\Big\{\kappa\ge1:\ \sup_{t\in\R}\|q(t)\|_{E^\eps_\kappa}\le\delta_0\Big\}.
\]
\end{definition}

Lemma~\ref{lem:kappa} below shows that $\kappa(q_0)<\infty$ for every Schwartz datum, that the infimum is attained, and that $\kappa(q_0)=1$ whenever $\Mass(q_0)\le\delta_0^2$. It also shows that, apart from the normalization $\kappa\ge1$, the frequency scale transforms linearly under the scaling \eqref{scaling}. Our first result makes \eqref{HKNV-LS} quantitative.

\begin{theorem}[Quantitative local smoothing]\label{thm:A}
There is a constant $C_1$, depending only on the fixed parameter $\eps$, such that for every $q_0\in\Sch$ and every $T>0$ the solution $q=\Phi(\cdot)q_0$ satisfies
\begin{equation}\label{thmA}
\LS_T(q)\le C_1\big(1+\kappa(q_0)^2\,T\big)\,\Mass(q_0).
\end{equation}
Conversely, the power of $\kappa$ in \eqref{thmA} cannot be lowered below one. More precisely, there are constants $c_*>0$, $\kappa_s>0$ and $\lambda_0\ge1$, and Schwartz data $q_{0,\lambda}$, $\lambda\ge1$, such that $\Mass(q_{0,\lambda})=2\pi$, $\kappa(q_{0,\lambda})=\max\{1,\kappa_s\lambda\}$, and
\begin{equation}\label{thmA-lower}
\LS_T\big(\Phi(\cdot)q_{0,\lambda}\big)\ge c_*\,\kappa(q_{0,\lambda})\,T\,\Mass(q_{0,\lambda})\qquad\text{for all }T>0\text{ and all }\lambda\ge\lambda_0 .
\end{equation}
\end{theorem}

Theorem~\ref{thm:A} recovers \cite[Theorem~1.5]{HKNV}. Indeed, if $Q\subset\Sch$ is bounded and equicontinuous, then \cite[Corollary~1.9]{HKNV} gives $\sup_{q_0\in Q}\kappa(q_0)<\infty$; see Remark~\ref{rem:constant}. The new information is twofold. All dependence on the data beyond the mass is carried by the single number $\kappa(q_0)$, and this dependence is polynomial. When the mass is small, $\kappa(q_0)=1$ and we obtain uniform estimates.

\begin{corollary}[Uniform local smoothing for small mass]\label{cor:small}
Let $m_0:=\delta_0^2$. For every $q_0\in L^2(\R)$ with $\Mass(q_0)\le m_0$ and every $T>0$, the solution $q=\Phi(\cdot)q_0$ of \cite{HKNV} satisfies
\begin{align}
\LS_T(q)&\le C_1(1+T)\,\Mass(q_0),\label{small1}\\
\sup_{h\in\R}\int_{-R^2}^{R^2}\big\|W\big((\cdot-h)/R\big)\,q(t)\big\|_{\Hdot}^2\,dt&\le 2C_1\,R\,\Mass(q_0)\qquad\text{for all }R>0,\label{small2}\\
\sup_{h\in\R}\big\|W_h\,q\big\|_{L^4([-T,T]\times\R)}^4&\le C_2(1+T)\,\Mass(q_0)^2 .\label{small3}
\end{align}
In particular $\Cm(m_0)\le 2C_1$ and $m^*\ge m_0>0$.
\end{corollary}

The estimate \eqref{small2} is scale invariant. It is the nonlinear counterpart, on parabolic cylinders, of the Kato smoothing estimate for the free Schr\"odinger group. The point of Corollary~\ref{cor:small} is not its difficulty, since it follows from the proof of \cite[Proposition~4.1]{HKNV}. Its point is that it locates the obstruction of \cite[Proposition~1.11]{HKNV} at positive mass. Our second result analyses this obstruction.

To state it, recall from \cite[Section~1.3]{HKNV} that for each $\vartheta\in(0,\pi/2)$ DNLS has a soliton of mass $8\vartheta$ travelling with velocity $-2\cot(2\vartheta)$. The soliton with $\vartheta=\pi/4$ is stationary. We denote by $\sol_\vartheta$ the family of all solutions obtained from it by translations in space and time, phase rotations and the scaling \eqref{scaling}; see \eqref{solfamily}.

\begin{theorem}[The stationary-soliton threshold]\label{thm:B}
There are absolute constants $c>0$, $C\ge1$ and $m_1\in[\pi,2\pi)$ such that the following hold.
\begin{enumerate}
\item[(i)] $\Cm(m)=\infty$ for every $m\ge2\pi$.
\item[(ii)] $\Cm(m)\ge c\,(2\pi-m)^{-1}$ for every $m\in[m_1,2\pi)$. In particular $\Cm(m)\to\infty$ as $m\uparrow2\pi$.
\item[(iii)] Every soliton $u\in\sol_\vartheta$ with $\vartheta\in(0,\pi/4)$ satisfies
\begin{equation}\label{thmB-iii}
\LS_1(u)\le C\Big(1+\frac{1}{2\pi-\Mass(u)}\Big)\Mass(u).
\end{equation}
\end{enumerate}
Consequently $m_0\le m^*\le2\pi$. Moreover, whatever the value of $m^*$, the constant $\Cm(m)$ is at least of order $(2\pi-m)^{-1}$ near $2\pi$, and the soliton family alone produces exactly this rate.
\end{theorem}

Part (i) is a reformulation of \cite[Proposition~1.11]{HKNV}. We include a proof because the same lower bound, applied to slowly moving solitons, gives (ii). Part (iii) shows that solitons alone cannot push $m^*$ below $2\pi$. Their contribution to local smoothing is controlled by the time they need to cross a window, which is of order $(\lambda\cot 2\vartheta)^{-1}$, against a local smoothing density of order $\lambda$; see Figure~\ref{fig:spacetime}. These considerations lead to the following conjecture, which is the main open problem of this paper.

\begin{conjecture}[Mass threshold for uniform local smoothing]\label{conj:A}
$m^*=2\pi$. More precisely, $\Cm(m)<\infty$ for every $m<2\pi$, and
\[
\Cm(m)\simeq(2\pi-m)^{-1}\qquad\text{as }m\uparrow2\pi .
\]
\end{conjecture}

In Section~\ref{sec:discussion} we explain why we believe this conjecture, using the trace formula for the Kaup--Newell scattering data. For Schwartz data of mass $m<2\pi$, every eigenvalue $\zeta_j$ of the Kaup--Newell problem satisfies $\arg\zeta_j\le m/4<\pi/2$. Consequently, no component of the discrete spectrum can be stationary; see Remark~\ref{rem:spectral}. We also formulate a sharp form of Theorem~\ref{thm:A} (Conjecture~\ref{conj:B}), in which $\kappa^2$ is replaced by $\kappa$, and a global-in-time question for small mass (Question~\ref{q:global}).

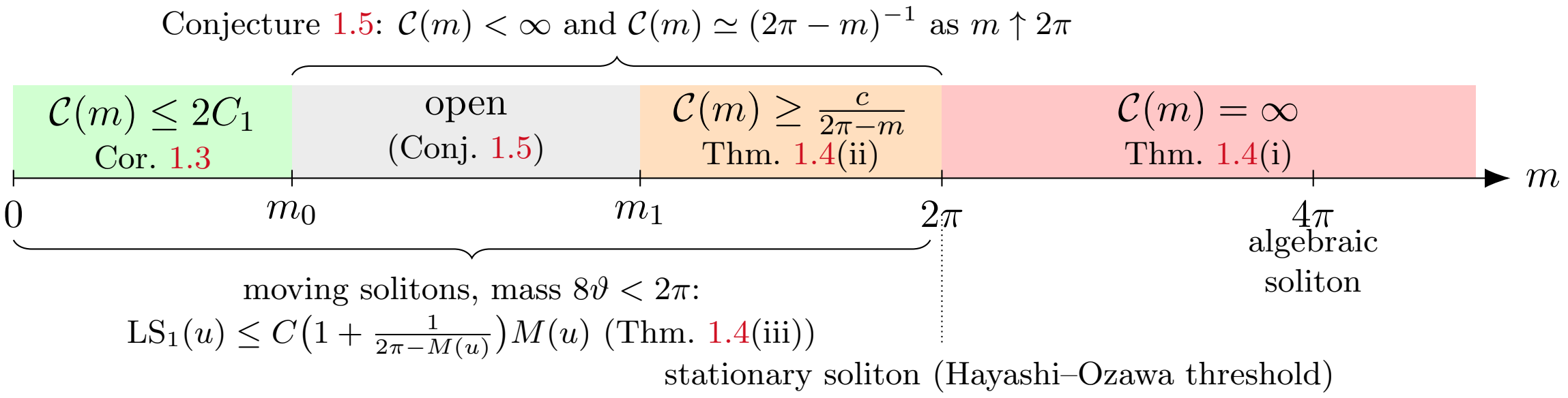
\begin{figure}[H]
\centering
\begin{tikzpicture}[x=2.25cm,y=1cm,font=\small]
\fill[green!18] (0,0) rectangle (1.2,0.9);
\fill[gray!15]  (1.2,0) rectangle (2.7,0.9);
\fill[orange!25] (2.7,0) rectangle (4.0,0.9);
\fill[red!22]   (4.0,0) rectangle (6.3,0.9);
\node[align=center] at (0.6,0.45) {$\Cm(m)\le 2C_1$\\[-1pt]{\scriptsize Cor.~\ref{cor:small}}};
\node[align=center] at (1.95,0.45) {open\\[-1pt]{\scriptsize (Conj.~\ref{conj:A})}};
\node[align=center] at (3.35,0.45) {$\Cm(m)\ge \frac{c}{2\pi-m}$\\[-1pt]{\scriptsize Thm.~\ref{thm:B}(ii)}};
\node[align=center] at (5.15,0.45) {$\Cm(m)=\infty$\\[-1pt]{\scriptsize Thm.~\ref{thm:B}(i)}};
\draw[-{Latex[length=2.5mm]}] (0,0) -- (6.45,0) node[right] {$m$};
\foreach \x/\l in {0/0,1.2/m_0,2.7/m_1,4.0/2\pi,5.6/4\pi}
  {\draw (\x,0.08) -- (\x,-0.08) node[below] {$\l$};}
\draw[decorate,decoration={brace,amplitude=5pt}] (1.2,1.0) -- (4.0,1.0)
  node[midway,above=5pt,align=center] {\scriptsize Conjecture~\ref{conj:A}: $\Cm(m)<\infty$ and $\Cm(m)\simeq(2\pi-m)^{-1}$ as $m\uparrow2\pi$};
\draw[decorate,decoration={brace,mirror,amplitude=5pt}] (0,-0.6) -- (3.95,-0.6)
  node[midway,below=6pt,align=center] {\scriptsize moving solitons, mass $8\vartheta<2\pi$:\\[-1pt] \scriptsize $\LS_1(u)\le C\big(1+\frac{1}{2\pi-\Mass(u)}\big)\Mass(u)$ (Thm.~\ref{thm:B}(iii))};
\node[align=center,below] at (4.25,-1.62) {\scriptsize stationary soliton (Hayashi--Ozawa threshold)};
\draw[densely dotted] (4.0,-0.35) -- (4.0,-1.62);
\node[align=center,below] at (5.6,-0.35) {\scriptsize algebraic\\[-2pt]\scriptsize soliton};
\end{tikzpicture}
\caption{The uniform local smoothing constant $\Cm(m)$ along the mass axis (not to scale). Green: proved finite. Orange: proved lower bound, finiteness conjectured. Red: proved infinite. The soliton bound of Theorem~\ref{thm:B}(iii) matches the lower bound of Theorem~\ref{thm:B}(ii).}
\label{fig:massaxis}
\end{figure}

Figure~\ref{fig:massaxis} summarizes what is proved and what is conjectured about $\Cm(m)$, and Table~\ref{tab:LS} places our results among the known local smoothing statements.

\begin{table}[H]
\centering
\small
\renewcommand{\arraystretch}{1.25}
\begin{tabular}{@{}>{\raggedright}p{3.6cm}>{\raggedright}p{5.4cm}>{\raggedright}p{2.6cm}p{2.3cm}@{}}
\toprule
Setting & Statement & Constant & Source\\
\midrule
Free Schr\"odinger group & $\sup_x\int_\R|D_x^{1/2}e^{it\partial_x^2}u_0(x)|^2dt=c\|u_0\|_{L^2}^2$ & exact & \cite{Kato83,KPV91}\\
DNLS, $Q\subset\Sch$ bounded and equicontinuous & $\LS_T(q)\lesssim_{T,Q}\Mass(q_0)$ & implicit in $Q$ & \cite[Thm.~1.5]{HKNV}\\
DNLS, rescaled stationary solitons, $\Mass\equiv2\pi$ & local smoothing fails on balls & no bound & \cite[Prop.~1.11]{HKNV}\\
DNLS, $q_0\in\Sch$ & $\LS_T(q)\le C_1(1+\kappa(q_0)^2T)\Mass(q_0)$; power $\ge1$ necessary & explicit in $\kappa(q_0)$ & Thm.~\ref{thm:A} (new)\\
DNLS, $\Mass(q_0)\le m_0$ & uniform bound, Morrey form, localized $L^4_{t,x}$ & uniform & Cor.~\ref{cor:small} (implicit in \cite{HKNV})\\
DNLS, $m\uparrow2\pi$ & $\Cm(m)\ge c(2\pi-m)^{-1}$ & blows up & Thm.~\ref{thm:B}(i),(ii) (new)\\
Solitons of mass $<2\pi$ & $\LS_1(u)\le C\big(1+\frac{1}{2\pi-\Mass(u)}\big)\Mass(u)$ & sharp rate & Thm.~\ref{thm:B}(iii) (new)\\
DNLS, $m<2\pi$ & $\Cm(m)<\infty$, $\Cm(m)\simeq(2\pi-m)^{-1}$ & conjectured & Conj.~\ref{conj:A}, \ref{conj:spectral}\\
\bottomrule
\end{tabular}
\caption{Local smoothing statements on $\R$ and the position of the present results.}
\label{tab:LS}
\end{table}

\subsection{Related literature}\label{subsec:lit}

\emph{Local and global well-posedness.} Local well-posedness of \eqref{DNLS} in $H^s(\R)$ was proved for $s>3/2$ by Tsutsumi and Fukuda \cite{TsutsumiFukuda80}, and in $H^1$ by Hayashi and Ozawa \cite{Hayashi93,HayashiOzawa92,HayashiOzawa94} by means of a gauge transformation; see also Ozawa \cite{Ozawa96}. Takaoka \cite{Takaoka99} reached $H^s$, $s\ge1/2$, by combining the gauge transformation with the Fourier restriction method. Below $s=1/2$ the data-to-solution map fails to be $C^3$ \cite{Takaoka99,Takaoka01}, and it is not even uniformly continuous on bounded sets \cite{BiagioniLinares01}. Hence contraction arguments in $H^s$ cannot go further. Local results in Fourier--Lebesgue and modulation spaces, closer to the critical scaling, are due to Gr\"unrock \cite{Grunrock05} and to Guo, Ren and Wang \cite{GuoRenWang21}.

Global results first required a mass restriction. Using the gauged energy and the sharp Gagliardo--Nirenberg inequality, Hayashi and Ozawa \cite{HayashiOzawa94} and Ozawa \cite{Ozawa96} obtained global $H^1$ solutions for $\Mass<2\pi$. Colliander, Keel, Staffilani, Takaoka and Tao \cite{CKSTT01,CKSTT02} carried the $2\pi$ condition to $H^s$, $s>1/2$, by the $I$-method, and Miao, Wu and Xu \cite{MiaoWuXu11} reached $H^{1/2}$. Wu exploited the momentum to go beyond $2\pi$, first to $2\pi+\eps_*$ \cite{Wu13} and then to $4\pi$ \cite{Wu15} in $H^1$. Guo and Wu \cite{GuoWu17} obtained the same threshold $4\pi$ in $H^{1/2}$; see also Fukaya, Hayashi and Inui \cite{FukayaHayashiInui17}. The value $4\pi$ is the mass of the algebraic soliton \cite[(1.26)]{HKNV}, \cite{BahouriPerelman22}, for which all other polynomial conservation laws vanish; see the discussion following \cite[(1.26)]{HKNV}.

The inverse scattering approach, initiated by Lee \cite{Lee89}, removed the mass restriction for data in weighted spaces \cite{LiuPerrySulem16,JLPS18cpde,JLPS20anal,PelinovskyShimabukuro18,PSS17}. It also produced long-time asymptotics for soliton-free data \cite{LPS18aihp} and soliton resolution for generic weighted data \cite{JLPS18cmp}; see the survey \cite{JLPS20qam}. Bahouri and Perelman \cite{BahouriPerelman22} then proved global well-posedness in $H^{1/2}(\R)$ without restriction, with uniform-in-time bounds, by combining the integrable structure with profile decompositions. Bahouri, Leslie and Perelman \cite{BahouriLesliePerelman23} established uniform-in-time $H^s$ bounds for $s\ge1/2$.

At the level of low-regularity conservation laws, Klaus and Schippa \cite{KlausSchippa22} obtained a priori estimates in Besov spaces of positive regularity for small mass, following the approach of Killip, Vi\c{s}an and Zhang \cite{KillipVisanZhang18}. Tang and Xu \cite{TangXu21} derived the corresponding microscopic conservation laws for small mass and pointed out that their coercive currents should yield local smoothing in $H^s$. They did not state an $L^2$-level estimate. Killip, Ntekoume and Vi\c{s}an \cite{KNV23} proved equicontinuity of orbits for $\Mass<4\pi$ and global well-posedness in $H^{1/6}$ under the same restriction. Harrop-Griffiths, Killip and Vi\c{s}an \cite{HKVequi} then removed the mass restriction from the equicontinuity statement. The $L^2$ theorem \cite{HKNV} builds on all of these works, on the method of commuting flows of Killip and Vi\c{s}an \cite{KillipVisan19}, and on its second-generation form developed for the cubic NLS and mKdV \cite{HKVpi}. Table~\ref{tab:WP} summarizes this development.

\begin{table}[H]
\centering
\small
\renewcommand{\arraystretch}{1.2}
\begin{tabular}{@{}>{\raggedright}p{3.0cm}>{\raggedright}p{2.6cm}>{\raggedright}p{4.2cm}p{4.0cm}@{}}
\toprule
Space & Mass condition & Main tool & References\\
\midrule
$H^s$, $s>3/2$ (local) & none & energy method & \cite{TsutsumiFukuda80}\\
$H^1$ (local) & none & gauge transformation & \cite{Hayashi93,HayashiOzawa92,HayashiOzawa94}\\
$H^1$ (global) & $\Mass<2\pi$ & gauged energy, sharp Gagliardo--Nirenberg & \cite{HayashiOzawa94,Ozawa96}\\
$H^s$, $s\ge1/2$ (local) & none & gauge, Fourier restriction & \cite{Takaoka99}\\
$H^s$, $s>1/2$ (global) & $\Mass<2\pi$ & $I$-method & \cite{CKSTT01,CKSTT02}\\
$H^{1/2}$ (global) & $\Mass<2\pi$ & $I$-method, resonant decomposition & \cite{MiaoWuXu11}\\
$H^1$ (global) & $\Mass<2\pi+\eps_*$, then $\Mass<4\pi$ & momentum, variational arguments & \cite{Wu13,Wu15}\\
$H^{1/2}$ (global) & $\Mass<4\pi$ & $I$-method, momentum & \cite{GuoWu17}\\
weighted spaces (global) & none & inverse scattering & \cite{LiuPerrySulem16,JLPS18cpde,JLPS20anal,PSS17,PelinovskyShimabukuro18}\\
$H^{1/2}$ (global) & none & integrability, profile decompositions & \cite{BahouriPerelman22}\\
$H^{1/6}$ (global) & $\Mass<4\pi$ & conservation laws, equicontinuity & \cite{KNV23}\\
$L^2$ (global) & none & commuting flows, local smoothing & \cite{HKNV}\\
\bottomrule
\end{tabular}
\caption{Well-posedness of DNLS on the line (selection).}
\label{tab:WP}
\end{table}

\emph{The periodic problem.} On the circle, local well-posedness in $H^{1/2}$ is due to Herr \cite{Herr06}. Gr\"unrock and Herr \cite{GrunrockHerr08} treated Fourier--Lebesgue data, and Deng, Nahmod and Yue \cite{DengNahmodYue21} obtained the optimal Fourier--Lebesgue local theory. Global results under the restriction $\Mass<4\pi$ were obtained by Mosincat and Oh in $H^1(\mathbb{T})$ \cite{MosincatOh15} and by Mosincat in $H^{1/2}(\mathbb{T})$ \cite{Mosincat17}. Bahouri and Perelman first relaxed the restriction to $\Mass<8\pi$ \cite{BahouriPerelmanIMRN}, and very recently they proved global well-posedness in $H^1(\mathbb{T})$ with no mass restriction \cite{BahouriPerelmanT}. Well-posedness in $L^2(\mathbb{T})$ remains open. Local smoothing, which is the tool studied here, is not available on the circle.

\emph{Solitons.} The two-parameter soliton family and its orbital stability were studied by Colin and Ohta \cite{ColinOhta06}. Kwon and Wu \cite{KwonWu18} treated the algebraic soliton, and Le Coz and Wu \cite{LeCozWu18} multi-solitons. The coincidence of $2\pi$ and $4\pi$ with the masses of the static and algebraic solitons is noted in \cite{BahouriPerelman22}.

\emph{A parallel model.} The continuum Calogero--Moser model, also called the Calogero--Moser derivative NLS, is another integrable, mass-critical derivative Schr\"odinger equation. There the stationary soliton also has mass $2\pi$. G\'erard and Lenzmann \cite{GerardLenzmann} proved global well-posedness of smooth solutions with mass at most $2\pi$. Killip, Laurens and Vi\c{s}an \cite{KLV25} proved well-posedness in $L^2_+(\R)$ for focusing data of mass below $2\pi$. Above $2\pi$, Hogan and Kowalski \cite{HoganKowalski} constructed solutions with unbounded Sobolev norms and mass arbitrarily close to $2\pi$, and Kim, Kim and Kwon \cite{KimKimKwon} constructed smooth finite-time blow-up solutions with mass arbitrarily close to that of the soliton. Small-mass dispersive decay is proved in \cite{KLMV}. In that model, therefore, $2\pi$ is a \emph{well-posedness} threshold. For DNLS, well-posedness in $L^2$ holds for every mass \cite{HKNV}, and Conjecture~\ref{conj:A} proposes that $2\pi$ is instead the threshold for \emph{uniform smoothing}. A threshold given by a ground-state mass is also familiar from the focusing mass-critical NLS \cite{Weinstein83,Dodson15} and from the scattering theory below the ground state developed by Kenig and Merle \cite{KenigMerle06}. For recent one-dimensional and inhomogeneous instances of this philosophy, we refer to \cite{BaiMajdoubJia,MajdoubSaanouni}.

\emph{Local smoothing.} For the free Schr\"odinger equation, Kato-type smoothing goes back to Kato \cite{Kato83}, Constantin and Saut \cite{ConstantinSaut88}, Sj\"olin \cite{Sjolin87} and Vega \cite{Vega88}. In one dimension, Kenig, Ponce and Vega \cite{KPV91} proved the exact identity $\sup_x\int_\R|D_x^{1/2}e^{it\partial_x^2}u_0(x)|^2dt=c\|u_0\|_{L^2}^2$. The same authors used local smoothing to treat Schr\"odinger equations with derivative nonlinearities \cite{KPV93,KPV98}. The finite-scale behaviour of the smoothing quantity was studied by Vega and Visciglia \cite{VegaVisciglia07}. An exact identity and sharp constants at finite scale were obtained by Bai and the author \cite{BaiMajdoub}. The monotonicity approach to local smoothing used in \cite{HKNV} is discussed in \cite{Tao06}.

To the best of our knowledge, the literature contains no uniform local smoothing statement for DNLS on $L^2$-balls of positive mass beyond what is implicit in \cite[Proposition~4.1]{HKNV}. Likewise, the only negative statement we are aware of is \cite[Proposition~1.11]{HKNV}, which shows failure at mass $2\pi$; the remarks in \cite{TangXu21} concern small mass in $H^s$, $s>0$. The question whether $2\pi$ is sharp, and how the constant depends on the data, therefore appears to be open. We do not claim novelty for Corollary~\ref{cor:small}, which is implicit in the proof of \cite[Proposition~4.1]{HKNV}. Our contributions are the explicit dependence in Theorem~\ref{thm:A}, the two-sided soliton analysis of Theorem~\ref{thm:B}, which pins down the expected rate $(2\pi-m)^{-1}$, and the formulation of Conjectures~\ref{conj:A} and~\ref{conj:B} together with a spectral reformulation.

The paper is organized as follows. Section~\ref{sec:prelim} collects notation, the results of \cite{HKNV} that we use, and two elementary tools, namely a window-summation lemma for $\Hdot$ and bounds on the soliton profiles. Section~\ref{sec:quant} proves Theorem~\ref{thm:A} and Corollary~\ref{cor:small}, and Section~\ref{sec:solitons} proves Theorem~\ref{thm:B}. Section~\ref{sec:discussion} discusses the method, its limitations and the evidence for the conjectures. Section~\ref{sec:concl} contains concluding remarks and further directions.

\section{Preliminaries and technical tools}\label{sec:prelim}

\subsection{Notation}

We use the unitary Fourier transform $\hat f(\xi)=(2\pi)^{-1/2}\int_\R e^{-ix\xi}f(x)\,dx$, so that $\|\hat f\|_{L^2}=\|f\|_{L^2}$. Following \cite{HKNV}, we write
\[
\|f\|_{H^s}^2=\int_\R(4+\xi^2)^s|\hat f(\xi)|^2\,d\xi,\qquad
\|f\|_{\Hdot}^2=\int_\R|\xi|\,|\hat f(\xi)|^2\,d\xi .
\]
Since $|\xi|\le(4+\xi^2)^{1/2}\le2+|\xi|$, we have
\begin{equation}\label{H12split}
\|f\|_{\Hdot}^2\le\|f\|_{H^{1/2}}^2\le 2\|f\|_{L^2}^2+\|f\|_{\Hdot}^2 .
\end{equation}
We shall use the Gagliardo representation
\begin{equation}\label{Gagliardo}
\|f\|_{\Hdot}^2=\frac1{2\pi}\iint_{\R^2}\frac{|f(x)-f(y)|^2}{|x-y|^2}\,dx\,dy .
\end{equation}
It follows from Plancherel's theorem, applied for fixed $h=x-y$, and from the identity $\int_\R(2-2\cos h\xi)\,h^{-2}\,dh=2\pi|\xi|$. We also use the elementary interpolation bound
\begin{equation}\label{interp}
\|f\|_{H^{1/2}}^2\le\|f\|_{L^2}\big(2\|f\|_{L^2}+\|f'\|_{L^2}\big),
\end{equation}
which follows from the Cauchy--Schwarz inequality in $\xi$ and from $\|f\|_{H^1}^2=4\|f\|_{L^2}^2+\|f'\|_{L^2}^2\le(2\|f\|_{L^2}+\|f'\|_{L^2})^2$.

The weight $W=\sech^6(\cdot/99)$ from \eqref{weight} satisfies $0<W\le1$ and
\begin{equation}\label{Wprops}
|W'|\le\tfrac{6}{99}W\le\tfrac1{16}W,\qquad \frac{W(x)}{W(y)}\le 2^6e^{6|x-y|/99},\qquad \int_\R W^2\,dx=\frac{512}{7}<74 .
\end{equation}
The first inequality follows from $W'=-\frac{6}{99}\tanh(\cdot/99)\,W$. The second follows from $\sech a\le2e^{-|a|}$ and $\sech b\ge e^{-|b|}$. The value of $\int W^2=99\int_\R\sech^{12}$ is recorded in \cite[(2.1)]{HKNV}. Moreover, $W(s)\ge\sech^6(2/99)\ge\frac12$ for $|s|\le2$.

\subsection{Results from \texorpdfstring{\cite{HKNV}}{[HKNV]} used below}

We use \cite{HKNV} as a black box through three statements. The first is the global well-posedness theorem \cite[Theorem~1.1]{HKNV}. It provides \eqref{Phi} and the group property $\Phi(t+s)=\Phi(t)\Phi(s)$. Moreover, Schwartz data lead to global Schwartz solutions, which conserve the mass and the perturbation determinant $a(k;q)$.

The second is the equicontinuity theorem of \cite{HKVequi}, quoted as \cite[Theorem~1.7]{HKNV}. It states that if $Q\subset\Sch$ is bounded and equicontinuous in $L^2$, then so is the set
\begin{equation}\label{Qstar}
Q_*=\big\{q\in\Sch:\ a(k;q)\equiv a(k;\tilde q)\text{ for some }\tilde q\in Q\big\}
\end{equation}
of \cite[(1.11)]{HKNV}. We combine it with \cite[Lemma~2.1]{HKNV}, which states that $\sup_{q\in Q}\|q\|_{E^\eps_\kappa}\to0$ as $\kappa\to\infty$ for every bounded equicontinuous set $Q$. (This last fact is elementary: splitting the integral \eqref{Enorm} at a frequency $K$, we get $\|q\|_{E^\eps_\kappa}^2\le(K/2\kappa)^{2\eps}\|q\|_{L^2}^2+\int_{|\xi|>K}|\hat q|^2\,d\xi$, and equicontinuity makes the last term uniformly small.) For $q_0\in\Sch$ we also set
\[
\kappa_\flat(q_0):=\inf\Big\{\kappa>0:\ \sup_{t\in\R}\|\Phi(t)q_0\|_{E^\eps_\kappa}\le\delta_0\Big\}\in[0,\infty).
\]

\begin{lemma}\label{lem:kappa}
Let $q_0\in\Sch$ and $q(t)=\Phi(t)q_0$.
\begin{enumerate}
\item[(i)] $\kappa(q_0)<\infty$, and the infimum in Definition~\ref{def:kappa} is attained. More precisely, $\sup_t\|q(t)\|_{E^\eps_\kappa}\le\delta_0$ for every $\kappa\ge\kappa(q_0)$.
\item[(ii)] If $\Mass(q_0)\le\delta_0^2$, then $\kappa(q_0)=1$.
\item[(iii)] For every $\lambda>0$ and $\kappa>0$ one has $\|(q_0)_\lambda\|_{E^\eps_\kappa}=\|q_0\|_{E^\eps_{\kappa/\lambda}}$, where $(q_0)_\lambda=\lambda^{1/2}q_0(\lambda\,\cdot)$. Consequently, $\kappa(q_0)=\max\{1,\kappa_\flat(q_0)\}$ and $\kappa_\flat\big((q_0)_\lambda\big)=\lambda\,\kappa_\flat(q_0)$.
\end{enumerate}
\end{lemma}

\begin{proof}
For $\kappa>0$ set $F(\kappa):=\sup_{t\in\R}\|q(t)\|_{E^\eps_\kappa}$. The weight in \eqref{Enorm} is at most one and decreases in $\kappa$. Hence, by dominated convergence, each map $\kappa\mapsto\|q(t)\|_{E^\eps_\kappa}$ is continuous and nonincreasing. We infer that $F$ is nonincreasing and lower semicontinuous, being a supremum of continuous functions. Therefore the set $S:=\{\kappa>0:F(\kappa)\le\delta_0\}$ is closed in $(0,\infty)$ and is a half-line: either $S=(0,\infty)$, in which case $\kappa_\flat(q_0)=0$, or $S=[\kappa_\flat(q_0),\infty)$.

(i) The DNLS flow conserves $a(k;q)$. Hence the whole orbit $\{\Phi(t)q_0:t\in\R\}$ is contained in $Q_*$ with $Q=\{q_0\}$. A single Schwartz function forms a bounded equicontinuous set, so $Q_*$ is bounded and equicontinuous by \cite[Theorem~1.7]{HKNV}. By \cite[Lemma~2.1]{HKNV}, $F(\kappa)\to0$ as $\kappa\to\infty$, so $S$ is nonempty and $\kappa_\flat(q_0)<\infty$. Since $S$ is a closed half-line, $\kappa(q_0)=\inf(S\cap[1,\infty))$ belongs to $S$, and $S\supset[\kappa(q_0),\infty)$.

(ii) If $\Mass(q_0)\le\delta_0^2$, the weight in \eqref{Enorm} being at most one, we get $\|q(t)\|_{E^\eps_1}\le\|q(t)\|_{L^2}=\Mass(q_0)^{1/2}\le\delta_0$ for all $t$ by mass conservation. Thus $1\in S$ and $\kappa(q_0)=1$.

(iii) We have $\widehat{(q_0)_\lambda}(\xi)=\lambda^{-1/2}\hat q_0(\xi/\lambda)$, and the substitution $\xi=\lambda\eta$ gives
\[
\int\frac{|\xi|^{2\eps}\lambda^{-1}|\hat q_0(\xi/\lambda)|^2}{(4\kappa^2+\xi^2)^\eps}\,d\xi=\int\frac{|\eta|^{2\eps}|\hat q_0(\eta)|^2}{(4(\kappa/\lambda)^2+\eta^2)^\eps}\,d\eta .
\]
The identity $\kappa(q_0)=\max\{1,\kappa_\flat(q_0)\}$ follows from the half-line structure of $S$. Finally, by the scaling invariance \eqref{scaling} and uniqueness of Schwartz solutions, $\Phi(t)(q_0)_\lambda=\big(\Phi(\lambda^2t)q_0\big)_\lambda$. Taking the supremum in $t$ and using the identity just proved, we obtain $\sup_t\|\Phi(t)(q_0)_\lambda\|_{E^\eps_\kappa}=F(\kappa/\lambda)$, which yields $\kappa_\flat((q_0)_\lambda)=\lambda\kappa_\flat(q_0)$.
\end{proof}

The third statement is the local smoothing estimate \cite[Proposition~4.1, (4.1)]{HKNV}. It reads as follows. Let $Q\subset\Sch$ be bounded and equicontinuous in $L^2$, and assume that the set of DNLS orbits
\[
\big\{\Phi(t)q:\ |t|\le1,\ q\in Q\big\}
\]
is $\delta$-good for a sufficiently small $\delta$. (In \cite[Proposition~4.1]{HKNV} this set is written $\{e^{tJ\nabla H}q:|t|\le1,\ q\in Q\}$, where $H$ is the DNLS Hamiltonian, and it is also denoted by $Q_*$; it should not be confused with the isospectral set \eqref{Qstar} above.) Then $\|q\|_{X^{1/2}}\lesssim\|q(0)\|_{L^2}$ uniformly for $q(0)\in Q$. Here $\delta$-good means bounded, equicontinuous and subject to the smallness condition \cite[(1.19)]{HKNV}, namely $\sup\|q\|_{E^\eps_1}\le\delta$; see \cite[Definition~1.8]{HKNV}. The local smoothing norm is defined in \cite[(2.24), (2.25)]{HKNV} over the time interval $[-1,1]$, and by \cite[(2.32)]{HKNV} with $\kappa=1$ it satisfies
\[
\|q\|_{X^{1/2}}^2\simeq\sup_{h\in\R}\int_{-1}^1\|W_hq(t)\|_{H^{1/2}}^2\,dt ,
\]
with absolute implicit constants, because $\psi^{12}=W$ and $\|(4-\partial^2)^{-1/2}f\|_{H^{3/2}}=\|f\|_{H^{1/2}}$. Finally, a compact subset of $L^2$ is automatically bounded and equicontinuous, and the orbit $\{\Phi(t)q_0:|t-s|\le1\}$ of a Schwartz datum is compact. Applying \cite[Proposition~4.1]{HKNV} with $Q=\{\Phi(s)q_0\}$, we obtain the result in the following form.

\begin{proposition}[{\cite[Proposition~4.1]{HKNV}}]\label{prop:HKNV41}
There are constants $\delta_0\in(0,1]$ and $C_0\ge1$, depending only on $\eps$, such that the following holds. Let $q_0\in\Sch$ and $s\in\R$, and suppose that $\|\Phi(t)q_0\|_{E^\eps_1}\le\delta_0$ for all $|t-s|\le1$. Then
\begin{equation}\label{HKNV41}
\sup_{h\in\R}\int_{s-1}^{s+1}\big\|W_h\,\Phi(t)q_0\big\|_{H^{1/2}}^2\,dt\le C_0\,\Mass(q_0).
\end{equation}
\end{proposition}

Decreasing $\delta_0$ only strengthens the hypothesis, so the normalization $\delta_0\le1$ costs nothing. The time translation from $[-1,1]$ to $[s-1,s+1]$ uses only the group property and mass conservation: one applies the estimate to the Schwartz datum $\Phi(s)q_0$.

\begin{remark}\rm \label{rem:constant}
The essential point in Proposition~\ref{prop:HKNV41} is that $\delta_0$ and $C_0$ depend on nothing but $\eps$. We have checked this against the proof of \cite[Proposition~4.1]{HKNV} on pp.~895--896 of the published version. Combining the identity \cite[(4.4)]{HKNV} with \cite[Lemmas~4.2 and~4.9]{HKNV}, \cite[(2.4)]{HKNV} and \cite[(2.13)]{HKNV}, the authors arrive at \cite[(4.33)]{HKNV}, which for $\kappa=1$ reads
\[
\|q\|_{X^{1/2}}^2\lesssim\|q\|_{L^\infty_tE^\eps_1}^2\|q\|_{X^{1/2}}^2+\|q\|_{L^\infty_tE^\eps_1}^2\big[1+\|q\|_{L^\infty_tL^2_x}^2\big]+\|q\|_{L^\infty_tL^2_x}^2 .
\]
If $\|q\|_{L^\infty_tE^\eps_1}\le\delta$ with $\delta$ small compared with the implicit constant, the first term is absorbed into the left-hand side; this is legitimate since $\|q\|_{X^{1/2}}<\infty$ for Schwartz solutions. Using $\|q\|_{E^\eps_1}\le\min\{\delta,\|q\|_{L^2}\}$ and mass conservation, we deduce that $\|q\|_{X^{1/2}}^2\lesssim(2+\delta^2)\Mass(q_0)$. In the proof of (4.1), the set $Q$ enters only through the smallness \cite[(1.19)]{HKNV} along the trajectory and through \cite[Lemmas~4.2 and~4.9]{HKNV}, whose estimates are stated uniformly for $q$ in any set satisfying the hypotheses of \cite[Proposition~4.1]{HKNV}. Equicontinuity is used quantitatively in that section only through \cite[(2.3)]{HKNV}, and only for the refined estimate \cite[(4.2)]{HKNV}, which we do not need. This is consistent with the proof of \cite[Theorem~1.5]{HKNV} on p.~921, where the authors note that the resulting constant depends on $Q$ solely through the rescaling parameter. The paper \cite{HKNV} does not fix a normalization of the Fourier transform. A different normalization multiplies every norm in \eqref{Enorm}, in \cite[(1.19)]{HKNV} and in $H^{1/2}$ by the same absolute factor, so only the values of $\delta_0$ and $C_0$ change.

For a bounded equicontinuous set $Q\subset\Sch$, \cite[Corollary~1.9]{HKNV} provides $\lambda=\lambda(Q,\delta_0)$ such that the rescaled isospectral set $Q_*^\lambda=\{\sqrt\lambda\,q(\lambda\,\cdot): q\in Q_*\}$ is $\delta_0$-good. Since $Q_*$ contains every orbit issued from $Q$, Lemma~\ref{lem:kappa}(iii) gives $\kappa(q_0)\le\max\{1,1/\lambda\}$ for all $q_0\in Q$, which is the bound used after Theorem~\ref{thm:A}. In this sense Theorem~\ref{thm:A} is a quantitative form of the proof of \cite[Theorem~1.5]{HKNV}, in which a uniform rescaling makes all orbits $\delta$-good and the estimate is then summed over unit time intervals.
\end{remark}

\subsection{Summing small windows}

The following lemma converts local smoothing on windows of size $r\le1$ into local smoothing on unit windows. For $r>0$ and $\eta\in\R$ we set
\begin{equation}\label{phirEta}
\varphi_{r,\eta}(x):=W\big((x-\eta)/r\big).
\end{equation}

\begin{lemma}[Window summation]\label{lem:windows}
There is an absolute constant $C_3$ such that for every $r\in(0,1]$, every $h\in\R$ and every $f\in H^{1/2}(\R)$,
\begin{equation}\label{windows}
\|W_hf\|_{\Hdot}^2+\frac1r\|W_hf\|_{L^2}^2\le\frac{C_3}{r}\int_\R W(\eta-h)^2\Big(\|\varphi_{r,\eta}f\|_{\Hdot}^2+\frac1r\|\varphi_{r,\eta}f\|_{L^2}^2\Big)\,d\eta .
\end{equation}
\end{lemma}

\begin{proof}
By translation invariance we may take $h=0$. We first treat the $L^2$ term. Fix $x$ and let $|\eta-x|\le r$. Then $\varphi_{r,\eta}(x)\ge\frac12$, and \eqref{Wprops} gives $W(x)\le2^6e^{6/99}W(\eta)\le2^7W(\eta)$. We average the identity $|Wf(x)|^2=\frac1{2r}\int_{|\eta-x|\le r}|W(x)f(x)|^2\,d\eta$ over these $\eta$, insert the two bounds, and integrate in $x$. We obtain
\begin{equation}\label{L2part}
\|Wf\|_{L^2}^2\le\frac{2^{16}}{r}\int_\R W(\eta)^2\|\varphi_{r,\eta}f\|_{L^2}^2\,d\eta ,
\end{equation}
which controls the second term on the left of \eqref{windows}.

For the homogeneous term, we split the double integral in \eqref{Gagliardo} into the far region $|x-y|\ge r$ and the near region $|x-y|<r$. In the far region we use $|Wf(x)-Wf(y)|^2\le2|Wf(x)|^2+2|Wf(y)|^2$ and $\int_{|y-x|\ge r}|x-y|^{-2}dy=2/r$. It follows that the far region contributes at most $8r^{-1}\|Wf\|_{L^2}^2$, which is controlled by \eqref{L2part}.

In the near region we again average over the $\eta$ with $|\eta-x|\le r$. For such $\eta$, both $x$ and $y$ lie in $J_\eta=[\eta-2r,\eta+2r]$, where $\varphi_{r,\eta}\ge\frac12$. On $J_\eta$ we write $Wf=m_\eta g_\eta$ with $g_\eta=\varphi_{r,\eta}f$ and $m_\eta=W/\varphi_{r,\eta}$. From \eqref{Wprops} and $r\le1$ we get $|m_\eta|\le2^8W(\eta)$ on $J_\eta$. Moreover, $|\varphi_{r,\eta}'|\le\varphi_{r,\eta}/(16r)$ and $|W'|\le W/16$, so that
\[
|m_\eta'|\le\frac{|W'|}{\varphi_{r,\eta}}+\frac{W|\varphi_{r,\eta}'|}{\varphi_{r,\eta}^2}\le\frac{W}{4r}\le\frac{2^5W(\eta)}{r}\qquad\text{on }J_\eta .
\]
Writing $Wf(x)-Wf(y)=m_\eta(x)\big(g_\eta(x)-g_\eta(y)\big)+\big(m_\eta(x)-m_\eta(y)\big)g_\eta(y)$ and using the mean value theorem on the interval $J_\eta$, we deduce that
\[
|Wf(x)-Wf(y)|\le2^8W(\eta)|g_\eta(x)-g_\eta(y)|+2^5W(\eta)r^{-1}|x-y|\,|g_\eta(y)| .
\]
We square this inequality, divide by $|x-y|^2$ and integrate over the near region, using $\int_{|x-y|<r}dx=2r$ for the second term. We conclude that the near region contributes at most
\[
\frac1{2r}\int_\R W(\eta)^2\Big(2^{17}\cdot2\pi\,\|g_\eta\|_{\Hdot}^2+\frac{2^{12}}{r}\|g_\eta\|_{L^2}^2\Big)\,d\eta
\]
to the double integral. Dividing by $2\pi$ and collecting the three contributions gives \eqref{windows}.
\end{proof}

\subsection{The soliton family}

We now record what we use about the solitons. Following \cite[Section~1.3]{HKNV}, for $\vartheta\in(0,\pi/2)$ we put
\begin{equation}\label{q0theta}
q_0(x;\vartheta)=A(x;\vartheta)\,e^{-ix\cot(2\vartheta)},\qquad
A(x;\vartheta)=\sqrt{2\sin(2\vartheta)}\;\frac{Z^3}{|Z|^4},\qquad Z=\cosh(x-i\vartheta).
\end{equation}
Note that $Z=\cos\vartheta\cosh x-i\sin\vartheta\sinh x$, so that \eqref{q0theta} is exactly the formula of \cite{HKNV}. (We have also checked numerically that \eqref{soliton} below solves \eqref{DNLS} with this sign of the phase, and not with the opposite one, when $\vartheta\neq\pi/4$. Only $|\partial_x e^{-ix\cot 2\vartheta}|=|\cot2\vartheta|$ enters our estimates.) By \cite[(1.24)]{HKNV},
\begin{equation}\label{soliton}
q(t,x;\vartheta)=q_0\big(x+2\cot(2\vartheta)t;\vartheta\big)\,e^{it\csc^2(2\vartheta)}
\end{equation}
is a Schwartz solution of \eqref{DNLS} of mass $\|q_0(\cdot;\vartheta)\|_{L^2}^2=8\vartheta$ \cite[(1.25)]{HKNV}. For $\vartheta=\pi/4$ it reduces to the stationary soliton \cite[(1.28)]{HKNV}. We let $\sol_\vartheta$ be the set of all solutions
\begin{equation}\label{solfamily}
u(t,x)=e^{i\gamma}\lambda^{1/2}\,q\big(\lambda^2(t-t_0),\lambda(x-x_0);\vartheta\big),\qquad \lambda>0,\ \gamma,t_0,x_0\in\R .
\end{equation}
Every $u\in\sol_\vartheta$ has mass $8\vartheta$, and by \eqref{soliton}
\begin{equation}\label{solmod}
|u(t,x)|=\lambda^{1/2}\big|A\big(\lambda(x-x(t));\vartheta\big)\big|,\qquad |\partial_xu(t,x)|=\lambda^{3/2}\big|\partial_xq_0\big(\lambda(x-x(t));\vartheta\big)\big|,
\end{equation}
with centre $x(t)=x_0-2\lambda\cot(2\vartheta)(t-t_0)$. The second identity holds because the phase factor in \eqref{soliton} depends on $t$ only.

\begin{lemma}[Profile bounds]\label{lem:profile}
Let $\vartheta\in(0,\pi/4]$ and $b:=\cot(2\vartheta)\ge0$. Then for all $x\in\R$
\begin{equation}\label{profile}
|A(x;\vartheta)|^2=\frac{4\sin(2\vartheta)}{\cosh(2x)+\cos(2\vartheta)}\le 2\,(8\vartheta)\,e^{-2|x|},\qquad |\partial_xq_0(x;\vartheta)|\le(9+b)\,|A(x;\vartheta)| .
\end{equation}
Moreover, $\vartheta\mapsto q_0(\cdot;\vartheta)$ is continuous from $[\pi/8,\pi/4]$ into $H^1(\R)$.
\end{lemma}

\begin{proof}
Since $|\cosh(a+ib)|^2=\sinh^2a+\cos^2b$ and $|\sinh(a+ib)|^2=\sinh^2a+\sin^2b$, we have
\[
|Z|^2=\tfrac12\big(\cosh2x+\cos2\vartheta\big),\qquad |Z'|^2=\sinh^2x+\sin^2\vartheta\le\cosh^2x .
\]
As $|A|^2=2\sin(2\vartheta)|Z|^{-2}$, this gives the identity in \eqref{profile}. For $\vartheta\le\pi/4$ we have $\cos2\vartheta\ge0$. Hence $|A|^2\le4\sin(2\vartheta)/\cosh2x\le8\sin(2\vartheta)e^{-2|x|}$, and $\sin2\vartheta\le2\vartheta$ gives the first bound.

Next, $A'/A=3Z'/Z-2(|Z|^2)'/|Z|^2$. Because $\cosh2x\ge\cosh^2x$ and $\cos 2\vartheta\ge 0$, we have $|Z|^2\ge\frac12\cosh^2x$, and therefore $|Z'/Z|\le\sqrt2$. Furthermore, $(|Z|^2)'/|Z|^2=2\sinh2x/(\cosh2x+\cos2\vartheta)$ is bounded by $2$ in absolute value. We infer that $|A'|\le(3\sqrt2+4)|A|\le9|A|$, and the product rule applied to \eqref{q0theta} gives the second bound in \eqref{profile}.

For the continuity statement, note that $Z$ does not vanish on $\R\times[\pi/8,\pi/4]$. Hence $q_0(x;\vartheta)$ and $\partial_xq_0(x;\vartheta)$ are continuous in $(x,\vartheta)$ on this set. By \eqref{profile} (with $b\le1$), they are dominated by $C e^{-|x|}$ uniformly in $\vartheta$. Dominated convergence yields continuity in $H^1$.
\end{proof}

\section{Quantitative local smoothing}\label{sec:quant}

\subsection{Local smoothing at the natural scale}

The first step is to rescale Proposition~\ref{prop:HKNV41} to the scale on which the orbit is $\delta_0$-good. In the original variables, this gives local smoothing on windows of size $\kappa^{-1}$ over times of length $\kappa^{-2}$.

\begin{proposition}[Scale-adapted local smoothing]\label{prop:scaled}
Let $q_0\in\Sch$ and $\kappa\ge\kappa(q_0)$, $\kappa\ge1$. Then for every $s,\eta\in\R$,
\begin{equation}\label{scaled}
\int_{s-\kappa^{-2}}^{s+\kappa^{-2}}\big\|\varphi_{1/\kappa,\eta}\,q(t)\big\|_{\Hdot}^2\,dt\le\frac{C_0\,\Mass(q_0)}{\kappa},
\end{equation}
where $q=\Phi(\cdot)q_0$ and $\varphi_{1/\kappa,\eta}=W(\kappa(\cdot-\eta))$.
\end{proposition}

\begin{proof}
Let $v(\tau,y):=\kappa^{-1/2}q(\kappa^{-2}\tau,\kappa^{-1}y)$. This is the rescaling \eqref{scaling} with $\lambda=\kappa^{-1}$. Hence $v$ is the Schwartz solution with data $(q_0)_{1/\kappa}$, and $\Mass(v)=\Mass(q_0)$. By Lemma~\ref{lem:kappa}(iii) with $\lambda=1/\kappa$, and then by Lemma~\ref{lem:kappa}(i), we have $\|v(\tau)\|_{E^\eps_1}=\|q(\kappa^{-2}\tau)\|_{E^\eps_\kappa}\le\delta_0$ for all $\tau$. Proposition~\ref{prop:HKNV41} applied at time $\kappa^2s$, together with \eqref{H12split}, gives
\[
\int_{\kappa^2s-1}^{\kappa^2s+1}\big\|W_{\kappa\eta}\,v(\tau)\big\|_{\Hdot}^2\,d\tau\le C_0\,\Mass(q_0).
\]
Now $W_{\kappa\eta}(y)v(\tau,y)=\kappa^{-1/2}F_\tau(y/\kappa)$ with $F_\tau(x)=W(\kappa(x-\eta))\,q(\kappa^{-2}\tau,x)$. A change of variables in the Fourier integral shows that $\|\kappa^{-1/2}F(\cdot/\kappa)\|_{\Hdot}^2=\kappa^{-1}\|F\|_{\Hdot}^2$. Substituting $\tau=\kappa^2t$, we obtain $\kappa\int_{s-\kappa^{-2}}^{s+\kappa^{-2}}\|\varphi_{1/\kappa,\eta}q(t)\|_{\Hdot}^2\,dt\le C_0\Mass(q_0)$, which is \eqref{scaled}.
\end{proof}

\subsection{Proof of Theorem~\ref{thm:A}}

Let $q_0\in\Sch$, write $\kappa=\kappa(q_0)$ and $\Mass=\Mass(q_0)$, and fix $h\in\R$. Consider a time interval $B=[s-\kappa^{-2},s+\kappa^{-2}]$. We apply Lemma~\ref{lem:windows} with $r=\kappa^{-1}\le1$ to $f=q(t)$ for each $t\in B$, and we integrate over $B$. On the right-hand side, Proposition~\ref{prop:scaled} bounds the time integral of $\|\varphi_{1/\kappa,\eta}q(t)\|_{\Hdot}^2$ by $C_0\Mass/\kappa$ for every $\eta$. Since $\varphi_{1/\kappa,\eta}\le1$, mass conservation gives $\int_B\kappa\|\varphi_{1/\kappa,\eta}q(t)\|_{L^2}^2\,dt\le\kappa\cdot2\kappa^{-2}\Mass$. Hence, by Fubini's theorem and \eqref{Wprops},
\[
\int_B\|W_hq(t)\|_{\Hdot}^2\,dt\le C_3\,\kappa\int_\R W(\eta-h)^2\Big(\frac{C_0\Mass}{\kappa}+\frac{2\Mass}{\kappa}\Big)\,d\eta\le 74\,C_3(C_0+2)\,\Mass .
\]
Moreover, $\int_B\|W_hq(t)\|_{L^2}^2\,dt\le2\kappa^{-2}\Mass\le2\Mass$. Combining these two bounds with \eqref{H12split}, we get $\int_B\|W_hq(t)\|_{H^{1/2}}^2dt\le C_4\Mass$ with $C_4=74C_3(C_0+2)+4$. The interval $[-T,T]$ is covered by at most $1+\kappa^2T$ intervals of the form $B$. Summing their contributions and taking the supremum over $h$, we obtain \eqref{thmA} with $C_1=C_4$.

For the lower bound, let $q_s$ be the stationary soliton, that is \eqref{soliton} with $\vartheta=\pi/4$, and set $q_{0,\lambda}:=(q_s(0))_\lambda$. By \eqref{soliton}, $\Phi(t)q_{0,\lambda}=e^{i\lambda^2t}(q_s(0))_\lambda$, so that $\sup_t\|\Phi(t)q_{0,\lambda}\|_{E^\eps_\kappa}=\|q_s(0)\|_{E^\eps_{\kappa/\lambda}}$ by Lemma~\ref{lem:kappa}(iii). Set $\kappa_s:=\kappa_\flat(q_s(0))$. It is finite by Lemma~\ref{lem:kappa}(i). It is positive, because $\|q_s(0)\|_{E^\eps_\kappa}\to\|q_s(0)\|_{L^2}=\sqrt{2\pi}>\delta_0$ as $\kappa\downarrow0$, by monotone convergence. Lemma~\ref{lem:kappa}(iii) then gives $\kappa(q_{0,\lambda})=\max\{1,\kappa_s\lambda\}\le K_*\lambda$ for $\lambda\ge1$, where $K_*:=\max\{1,\kappa_s\}$. In particular $\kappa(q_{0,\lambda})$ grows linearly in $\lambda$.

On the other hand, $\Phi(\cdot)q_{0,\lambda}\in\sol_{\pi/4}$ with $x_0=t_0=0$, and its centre $x(t)\equiv0$ does not move. Proposition~\ref{prop:lower} below therefore gives $\|W\Phi(t)q_{0,\lambda}\|_{H^{1/2}}^2\ge c_5\lambda$ for all $t$ and all $\lambda\ge\lambda_0$. Integrating over $[-T,T]$, we find
\[
\LS_T\big(\Phi(\cdot)q_{0,\lambda}\big)\ge2c_5\lambda T\ge\frac{2c_5}{K_*}\,\kappa(q_{0,\lambda})\,T .
\]
Since $\Mass(q_{0,\lambda})=2\pi$, this is \eqref{thmA-lower} with $c_*=c_5/(\pi K_*)$. Finally, if \eqref{thmA} held with $\kappa^2$ replaced by $\kappa^\alpha$ for some $\alpha<1$, then $2c_5\lambda T\le C(1+K_*^\alpha\lambda^\alpha T)2\pi$ for all $\lambda\ge\lambda_0$, which is impossible as $\lambda\to\infty$. \qed

\begin{remark}\rm \label{rem:gap}
The exponent $2$ in \eqref{thmA} comes from parabolic scaling. Proposition~\ref{prop:scaled} controls windows of size $\kappa^{-1}$ only during times of order $\kappa^{-2}$. Moreover, the smallness that drives the bootstrap of \cite{HKNV} is available only at frequencies above $\kappa$, so one cannot simply enlarge the windows. The stationary soliton shows that the true growth is at least linear in $\kappa$; see Conjecture~\ref{conj:B}.
\end{remark}

\subsection{Proof of Corollary~\ref{cor:small}}

Let first $q_0\in\Sch$ with $\Mass(q_0)\le m_0=\delta_0^2$. By Lemma~\ref{lem:kappa}(ii), $\kappa(q_0)=1$, and \eqref{small1} is \eqref{thmA}.

For \eqref{small2}, let $R>0$ and let $q^{(R)}(t,x)=R^{1/2}q(R^2t,Rx)$ be the solution with data $(q_0)_R$, which has the same mass. Applying \eqref{small1} with $T=1$ to $q^{(R)}$ and using \eqref{H12split}, we obtain $\int_{-1}^1\|W_{h}q^{(R)}(t)\|_{\Hdot}^2dt\le2C_1\Mass(q_0)$ for every $h$. Now $W_h(x)q^{(R)}(t,x)=R^{1/2}G(Rx)$ with $G(y)=W((y-Rh)/R)\,q(R^2t,y)$, and $\|R^{1/2}G(R\,\cdot)\|_{\Hdot}^2=R\|G\|_{\Hdot}^2$. Hence the substitution $t'=R^2t$ turns the last bound into
\[
\int_{-R^2}^{R^2}\big\|W\big((\cdot-Rh)/R\big)q(t')\big\|_{\Hdot}^2\,dt'\le2C_1R\,\Mass(q_0).
\]
Since $Rh$ ranges over $\R$, this is \eqref{small2}.

For \eqref{small3}, we use the one-dimensional inequality $\|f\|_{L^4}\le C\|f\|_{\dot H^{1/4}}\le C\|f\|_{L^2}^{1/2}\|f\|_{\Hdot}^{1/2}$, which follows from Sobolev embedding and interpolation. We apply it to $f=W_hq(t)$. As $\|W_hq(t)\|_{L^2}^2\le\Mass(q_0)$, integration in time and \eqref{small1} give \eqref{small3} with $C_2=C^4C_1$.

Now let $q_0\in L^2$ with $0<\Mass(q_0)\le m_0$. Choose Schwartz functions $p_n\to q_0$ in $L^2$, with $p_n\neq0$, and set $q_{0,n}=\|q_0\|_{L^2}\,p_n/\|p_n\|_{L^2}$. Then $q_{0,n}\in\Sch$, $\Mass(q_{0,n})=\Mass(q_0)$, and $q_{0,n}\to q_0$ in $L^2$. By the joint continuity of $\Phi$ \cite[Theorem~1.1]{HKNV} and the compactness of $[-T,T]$, we get $q_n(t)=\Phi(t)q_{0,n}\to q(t)$ in $L^2$ uniformly for $t\in[-T,T]$. For fixed $h$ and $t$, it follows that $W_hq_n(t)\to W_hq(t)$ in $L^2$.

The norms $\|\cdot\|_{H^{1/2}}$ and $\|\cdot\|_{\Hdot}$ (with value $+\infty$ allowed) are lower semicontinuous with respect to $L^2$ convergence. Indeed, if $f_n\to f$ in $L^2$, choose a subsequence along which $\|f_n\|_{H^{1/2}}$ tends to $\liminf_n\|f_n\|_{H^{1/2}}$, and a further subsequence along which $\hat f_n\to\hat f$ almost everywhere. Fatou's lemma in $\xi$ then gives $\|f\|_{H^{1/2}}\le\liminf_n\|f_n\|_{H^{1/2}}$, and the same argument applies to $\Hdot$. Fatou's lemma in $t$ therefore transfers \eqref{small1} and \eqref{small2} from $q_n$ to $q$. Finally, $W_hq_n\to W_hq$ in $L^2_{t,x}([-T,T]\times\R)$. Hence a subsequence converges almost everywhere, and Fatou's lemma also transfers \eqref{small3}. The bound $\Cm(m_0)\le2C_1$ is \eqref{small1} with $T=1$. \qed

\section{Solitons and the threshold \texorpdfstring{$2\pi$}{2pi}}\label{sec:solitons}

\begin{figure}[H]
\centering
\begin{tikzpicture}[x=0.62cm,y=2.2cm,font=\small]
\fill[blue!8] (-2.6,-1) rectangle (2.6,1);
\draw[blue!50,densely dashed] (-2.6,-1) -- (-2.6,1);
\draw[blue!50,densely dashed] (2.6,-1) -- (2.6,1);
\node[blue!60!black,align=center] at (0,1.17) {\scriptsize window $W_h$ (width $\sim99$)};
\draw[-{Latex[length=2.2mm]}] (-8.5,-1) -- (8.8,-1) node[right] {$x$};
\draw[-{Latex[length=2.2mm]}] (-8.5,-1) -- (-8.5,1.25) node[above] {$t$};
\node[left] at (-8.5,-1) {$-1$};
\node[left] at (-8.5,1) {$1$};
\draw[gray!60] (-8.5,1) -- (8.5,1);
\draw[red!75!black,very thick] (0,-1) -- (0,1);
\node[red!75!black,align=left,right] at (0.1,-0.82) {\scriptsize stationary ($b=0$)};
\draw[orange!85!black,very thick] (1.5,-1) -- (-1.5,1);
\node[orange!85!black,align=left,left] at (-1.6,0.88) {\scriptsize slow: $\lambda b\lesssim1$};
\draw[green!50!black,very thick] (8.3,-0.42) -- (-8.3,0.42);
\node[green!50!black,align=left] at (6.4,-0.55) {\scriptsize fast: $\lambda b\gg1$};
\draw[green!50!black,decorate,decoration={brace,amplitude=4pt}] (3.0,0.13) -- (3.0,-0.13);
\node[green!50!black,right] at (3.15,0.0) {\scriptsize $\sim(\lambda b)^{-1}$};
\end{tikzpicture}
\caption{Trajectories $x(t)=x_0-2\lambda b\,t$, $b=\cot(2\vartheta)$, of rescaled solitons in the space-time window $[-1,1]\times\operatorname{supp}W_h$ (schematic). Each soliton carries a local smoothing density of order $\lambda$, so its contribution to $\LS_1$ is of order $\lambda\cdot\min\{1,(\lambda b)^{-1}\}\le b^{-1}$. This gives Proposition~\ref{prop:upper}. The slow solitons with $\lambda b\simeq1$ and $\lambda\simeq b^{-1}$ attain this order, which gives Theorem~\ref{thm:B}(ii).}
\label{fig:spacetime}
\end{figure}
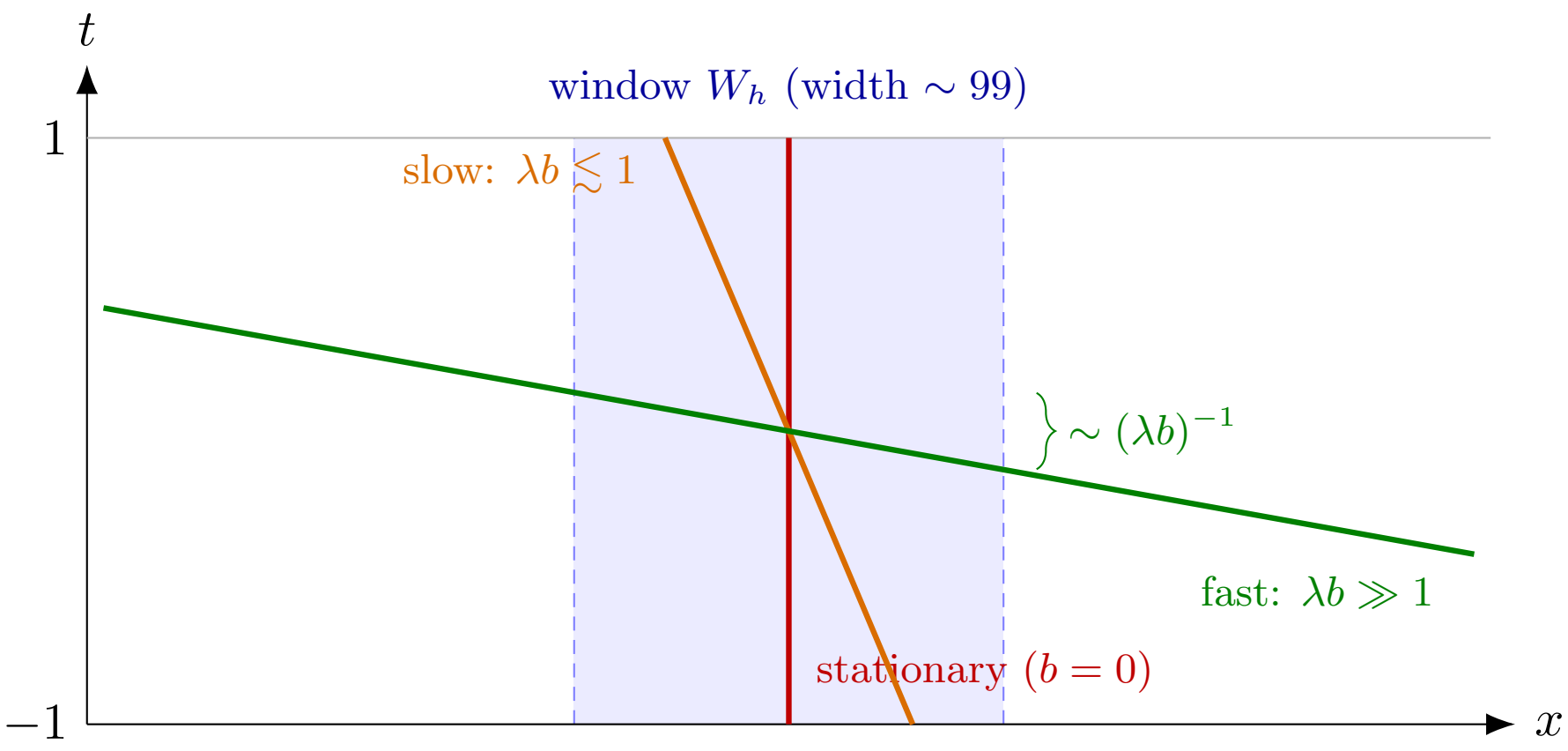

The two propositions of this section quantify a simple geometric mechanism, illustrated in Figure~\ref{fig:spacetime}. A soliton of scale $\lambda$ carries a local smoothing density of order $\lambda$, and it contributes only while it stays in the window. A soliton of mass $8\vartheta<2\pi$ moves with speed $2\lambda\cot(2\vartheta)$, so it crosses a unit window in a time of order $(\lambda\cot2\vartheta)^{-1}$. Its total contribution is therefore of order $(\cot2\vartheta)^{-1}$, uniformly in $\lambda$. By contrast, the stationary soliton never leaves the window.

\subsection{An upper bound for moving solitons}

\begin{proposition}\label{prop:upper}
Let $\vartheta\in(0,\pi/4)$, $b=\cot(2\vartheta)>0$ and $u\in\sol_\vartheta$. Then
\begin{equation}\label{upper}
\LS_1(u)\le\Big(86+\frac{666}{b}\Big)\Mass(u).
\end{equation}
\end{proposition}

\begin{proof}
Write $\Mass=\Mass(u)=8\vartheta$, fix $h\in\R$, and let $\lambda$ and $x(t)$ be as in \eqref{solmod}. Define
\[
I(t):=\int_\R W_h(x)^2\,\lambda\,e^{-2\lambda|x-x(t)|}\,dx .
\]
Since $W_h\le1$, we have $I(t)\le1$. By \eqref{solmod} and \eqref{profile},
\[
\|W_hu(t)\|_{L^2}^2=\int W_h^2\,\lambda\,|A(\lambda(x-x(t)))|^2\,dx\le2\Mass\,I(t).
\]
Using $|W_h'|\le W_h/16$ and the second bound in \eqref{profile}, we also get
\begin{align*}
\|(W_hu(t))'\|_{L^2}&\le\tfrac1{16}\|W_hu(t)\|_{L^2}+\lambda(9+b)\Big(\int W_h^2\,\lambda\,|A(\lambda(x-x(t)))|^2dx\Big)^{1/2}\\
&\le\Big(\tfrac1{16}+\lambda(9+b)\Big)\sqrt{2\Mass\,I(t)} .
\end{align*}
By \eqref{interp}, it follows that $\|W_hu(t)\|_{H^{1/2}}^2\le2\Mass\,I(t)\big(3+\lambda(9+b)\big)$.

We now integrate in time. The centre $x(t)$ moves with speed $2\lambda b$. Denote by $t_x$ the time at which $x(t)=x$, so that $|x-x(t)|=2\lambda b|t-t_x|$. Fubini's theorem gives
\[
\int_{-1}^1I(t)\,dt\le\int_\R W_h(x)^2\Big(\int_\R\lambda e^{-4\lambda^2b|t-t_x|}\,dt\Big)dx=\frac{\|W\|_{L^2}^2}{2\lambda b}\le\frac{37}{\lambda b}.
\]
Combining this with $I\le1$, we obtain
\[
\int_{-1}^1\|W_hu(t)\|_{H^{1/2}}^2dt\le2\Mass\Big((3+9\lambda)\min\Big\{2,\frac{37}{\lambda b}\Big\}+\lambda b\cdot\frac{37}{\lambda b}\Big)\le2\Mass\Big(6+\frac{333}{b}+37\Big),
\]
where we used the bound $3\cdot2$ for the first summand and $9\lambda\cdot37/(\lambda b)$ for the second. Taking the supremum over $h$ gives \eqref{upper}.
\end{proof}

The bound \eqref{upper} is uniform in the scaling parameter $\lambda$ and in the translation and phase parameters. It degenerates only as $b\downarrow0$, that is, as $\vartheta\uparrow\pi/4$. This is the regime of slowly moving solitons with mass close to $2\pi$. For small mass the bound is uniform as well, because small solitons travel fast.

\subsection{A lower bound for slow solitons}

\begin{proposition}\label{prop:lower}
There are absolute constants $c_5>0$ and $\lambda_0\ge1$ such that the following holds. Let $\vartheta\in[\pi/8,\pi/4]$ and $\lambda\ge\lambda_0$, and let $u\in\sol_\vartheta$ have parameters $x_0=t_0=0$ in \eqref{solfamily}. Then for every $t$ with $|x(t)|\le2$,
\begin{equation}\label{lower}
\|W\,u(t)\|_{H^{1/2}}^2\ge c_5\,\lambda .
\end{equation}
\end{proposition}

\begin{proof}
By \eqref{solfamily} and \eqref{soliton}, $W(x)u(t,x)=\omega\,\lambda^{1/2}G(\lambda(x-x(t)))$ for a unimodular constant $\omega$, where $G(y)=W(x(t)+y/\lambda)\,q_0(y;\vartheta)$. By \eqref{H12split}, translation invariance and the scaling of $\Hdot$,
\[
\|Wu(t)\|_{H^{1/2}}^2\ge\|Wu(t)\|_{\Hdot}^2=\lambda\,\|G\|_{\Hdot}^2 .
\]
Set $\mu:=\min_{\vartheta\in[\pi/8,\pi/4]}\|q_0(\cdot;\vartheta)\|_{\Hdot}$. This minimum is attained, by Lemma~\ref{lem:profile} and the continuity of the $\Hdot$ norm on $H^1$. It is positive, because a function with vanishing $\Hdot$ norm vanishes, while $\|q_0(\cdot;\vartheta)\|_{L^2}^2=8\vartheta>0$. Since $|x(t)|\le2$, we have $W(x(t))\ge\frac12$, and hence $\|W(x(t))q_0\|_{\Hdot}\ge\mu/2$.

It remains to estimate the difference $D(y)=\big(W(x(t)+y/\lambda)-W(x(t))\big)q_0(y;\vartheta)$. By \eqref{Wprops}, $|W(x(t)+y/\lambda)-W(x(t))|\le|y|/(16\lambda)$. By Lemma~\ref{lem:profile} with $b\le1$ and $16\vartheta\le4\pi$, we also have $|q_0|+|\partial_yq_0|\le 12\sqrt{4\pi}\,e^{-|y|}$. Moreover, $|D'(y)|\le\lambda^{-1}|W'||q_0|+(16\lambda)^{-1}|y|\,|\partial_yq_0|$. We deduce that
\[
\|D\|_{\Hdot}^2\le\|D\|_{L^2}^2+\|D'\|_{L^2}^2\le\frac{C_6^2}{\lambda^2},
\]
with an absolute constant $C_6$. We choose $\lambda_0:=\max\{1,4C_6/\mu\}$. Then $\|G\|_{\Hdot}\ge\mu/2-C_6/\lambda\ge\mu/4$, and \eqref{lower} holds with $c_5=\mu^2/16$.
\end{proof}

\subsection{Proof of Theorem~\ref{thm:B}}

(i) Let $m\ge2\pi$. For $\lambda\ge\lambda_0$, let $u_\lambda\in\sol_{\pi/4}$ be the rescaled stationary soliton with $x_0=t_0=0$. It has mass $2\pi\le m$ and centre $x(t)\equiv0$. By Proposition~\ref{prop:lower}, $\LS_1(u_\lambda)\ge2c_5\lambda$. Since $u_\lambda=\Phi(\cdot)u_\lambda(0)$, we get $\Cm(m)\ge c_5\lambda/\pi$ for every $\lambda\ge\lambda_0$, and therefore $\Cm(m)=\infty$.

(ii) Let $\lambda_0$ be as in Proposition~\ref{prop:lower} and set $m_1:=\max\{\pi,\ 2\pi-4\arctan(1/\lambda_0)\}$. Given $m\in[m_1,2\pi)$, let $\vartheta=m/8\in[\pi/8,\pi/4)$, so that $b=\cot(2\vartheta)=\tan\big((2\pi-m)/4\big)\in(0,1/\lambda_0]$. Take $u\in\sol_\vartheta$ with $x_0=t_0=0$ and $\lambda=1/b\ge\lambda_0$. Then $|x(t)|=2\lambda b|t|\le2$ for $|t|\le1$, and Proposition~\ref{prop:lower} gives $\LS_1(u)\ge2c_5/b$. Since $(2\pi-m)/4\le\pi/4$ and $\tan s\le4s/\pi$ on $[0,\pi/4]$ by convexity, we have $b\le(2\pi-m)/\pi$. As $\Mass(u)=m$, we conclude that
\[
\Cm(m)\ge\frac{\LS_1(u)}{\Mass(u)}\ge\frac{2c_5}{m\,b}\ge\frac{2c_5\pi}{2\pi(2\pi-m)}=\frac{c_5}{2\pi-m}.
\]

(iii) For $u\in\sol_\vartheta$ with $\vartheta\in(0,\pi/4)$ we have $\Mass(u)=8\vartheta$ and $b=\tan\big((2\pi-\Mass(u))/4\big)\ge(2\pi-\Mass(u))/4$. Hence $666/b\le2664/(2\pi-\Mass(u))$, and Proposition~\ref{prop:upper} gives \eqref{thmB-iii} with $C=2664$.

The final assertions follow from Corollary~\ref{cor:small} and parts (i)--(iii), with $c=c_5$. \qed

\begin{remark}\rm \label{rem:fast}
Proposition~\ref{prop:upper} also covers solitons of very small mass. As $\vartheta\downarrow0$, the amplitude of $q_0(\cdot;\vartheta)$ is of order $\vartheta^{1/2}$, while its frequency $\cot(2\vartheta)$ is of order $\vartheta^{-1}$. These solitons are thus high-frequency wave packets of small mass, travelling at a speed comparable to their frequency. They are invisible to the threshold question. They do show, however, that small-mass solutions need not scatter. This is why the global-in-time question below is phrased in terms of smoothing and not of scattering; compare with the weighted-space asymptotics of \cite{HayashiOzawaMA94,LPS18aihp}.
\end{remark}

\section{Discussion: methods, obstacles and conjectures}\label{sec:discussion}

\subsection{What the method gives and where it stops}

The proof of Theorem~\ref{thm:A} uses \cite{HKNV} in an essentially optimal way. Smallness of the $E^\eps$ norm at frequency scale $\kappa$ is the only source of smallness in their bootstrap, and rescaling converts it into local smoothing at the natural parabolic scale. All further losses come from geometry, namely from summing $\kappa$ small windows into a unit window and $\kappa^2T$ short time intervals into $[-T,T]$. Lemma~\ref{lem:windows} shows that the first step costs nothing. The second step is where the factor $\kappa^2$ appears.

For data of mass below $m_0$ the frequency scale is irrelevant, and Corollary~\ref{cor:small} follows. Between $m_0$ and $2\pi$ the method gives no uniform information. The reason is that the $E^\eps_1$ norm of a datum of moderate mass need not be small, and making it small requires a rescaling that depends on the frequency distribution of the datum. This is not a technical defect. By Theorem~\ref{thm:B}(i), the same mechanism must fail at mass $2\pi$. Any proof of Conjecture~\ref{conj:A} must therefore use a property of the region $\Mass<2\pi$ that the stationary soliton violates.

\subsection{Heuristic evidence for Conjecture~\ref{conj:A}}

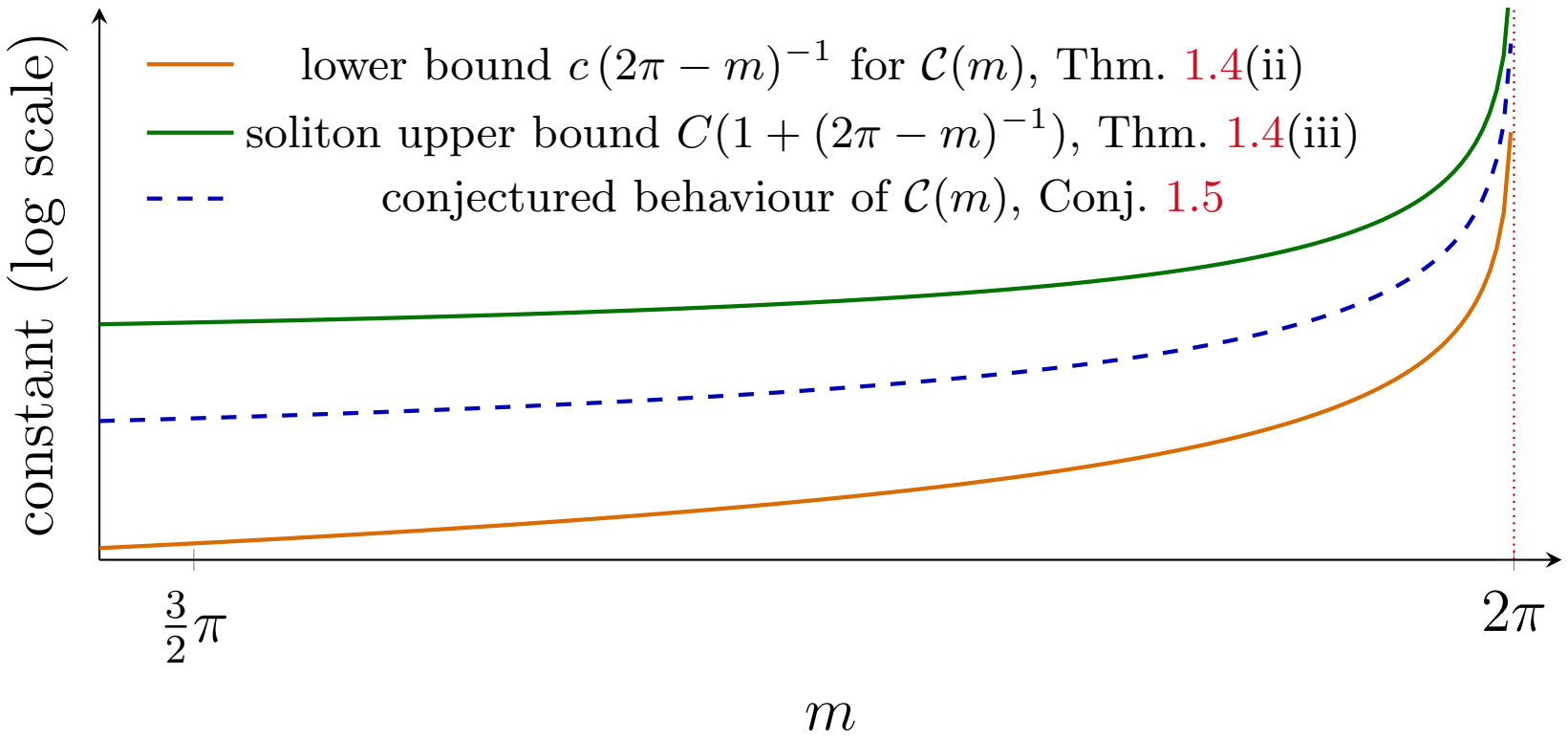
\begin{figure}[H]
\centering
\begin{tikzpicture}
\begin{axis}[
  width=0.72\textwidth, height=5.4cm,
  xmin=4.6, xmax=6.34, ymin=0.3, ymax=1000,
  ymode=log, axis lines=left,
  xlabel={$m$}, ylabel={constant (log scale)},
  xtick={4.712389,6.283185}, xticklabels={$\tfrac32\pi$,$2\pi$},
  ytick=\empty,
  legend style={font=\scriptsize,at={(0.02,0.98)},anchor=north west,draw=none,fill=none},
  every axis plot/.append style={thick,samples=200}]
\addplot[orange!85!black,domain=4.6:6.2795] {0.6/(6.283185-x)};
\addlegendentry{lower bound $c\,(2\pi-m)^{-1}$ for $\Cm(m)$, Thm.~\ref{thm:B}(ii)}
\addplot[green!45!black,domain=4.6:6.2795] {6*(1+1/(6.283185-x))};
\addlegendentry{soliton upper bound $C(1+(2\pi-m)^{-1})$, Thm.~\ref{thm:B}(iii)}
\addplot[blue!70!black,dashed,domain=4.6:6.2795] {2.2/(6.283185-x)+1};
\addlegendentry{conjectured behaviour of $\Cm(m)$, Conj.~\ref{conj:A}}
\draw[red!70!black,densely dotted] (axis cs:6.283185,0.3) -- (axis cs:6.283185,1000);
\end{axis}
\end{tikzpicture}
\caption{Schematic behaviour near $m=2\pi$ (constants not to scale). Both proved bounds blow up at the same rate $(2\pi-m)^{-1}$. Conjecture~\ref{conj:A} asserts that no other mechanism produces a faster rate below $2\pi$.}
\label{fig:rate}
\end{figure}

A unit window can accumulate a large amount of local smoothing only if a solution concentrated at a high frequency $\lambda$ remains near the window for a time much longer than $\lambda^{-1}$. Linear wave packets at frequency $\lambda$ move with speed of order $\lambda$ and do not do this, but coherent structures could. Proposition~\ref{prop:upper} shows that single solitons of mass $8\vartheta<2\pi$ cannot. Indeed, their velocity $2\lambda\cot(2\vartheta)$ is proportional to their frequency, with a coefficient bounded below in terms of $2\pi-8\vartheta$. Theorem~\ref{thm:B}(ii) shows that this coefficient is the only source of growth among solitons. Figure~\ref{fig:rate} illustrates the resulting two-sided picture near $2\pi$.

The scattering theory of the Kaup--Newell problem gives a more structural version of this observation. For $u\in\Sch$, let $\tilde a_u(\zeta)$ be the transmission coefficient in the variable $\zeta=k^2$, where $k$ is the Kaup--Newell spectral parameter, and let $\zeta_j$ be its zeros in the upper half-plane, with $\arg\zeta_j\in(0,\pi)$. Bahouri and Perelman record the trace formula \cite[Lemma~2.3(ii)]{BahouriPerelmanT}; see also \cite[(2.15)]{BahouriPerelman22}:
\begin{equation}\label{trace}
\Mass(u)=4\sum_j\arg\zeta_j+\frac2\pi\int_\R d\mu_u(\xi),
\end{equation}
where $\mu_u$ is a nonnegative measure accounting for the radiation.

\begin{remark}\rm [Spectral picture below $2\pi$]\label{rem:spectral}
By \eqref{trace} and the nonnegativity of $\mu_u$, every eigenvalue of a Schwartz function $u$ with $\Mass(u)<2\pi$ satisfies $0<\arg\zeta_j\le\Mass(u)/4<\pi/2$. For a single soliton the picture is explicit. In the normalization of \cite{BahouriPerelman22}, the bright soliton with parameters $(E,c)$, $E>0$, travels with velocity $c$. Its transmission coefficient has the single zero $\zeta_{E,c}=-c/4+i\sqrt E/2$, and its mass is given in \cite[(1.12)]{BahouriPerelman22}. A direct computation shows that $\cos(\arg\zeta_{E,c})=-c/\sqrt{c^2+4E}=\cos\big(\Mass(u_{E,c})/4\big)$, hence $\Mass(u_{E,c})=4\arg\zeta_{E,c}$. Consequently, there is no radiation in \eqref{trace}, the velocity $c=-4\operatorname{Re}\zeta_{E,c}$ is proportional to $\operatorname{Re}\zeta_{E,c}$, and the soliton is stationary if and only if its eigenvalue lies on the imaginary axis, that is, if and only if its mass equals $2\pi$. In the parametrization \eqref{q0theta} this corresponds to $\arg\zeta=2\vartheta$. The equation of \cite{BahouriPerelman22} is related to \eqref{DNLS} by the reflection $x\mapsto-x$, which reverses the direction of motion but does not affect these statements. For several solitons, breathers and general data, the analogous statement, namely that each coherent component attached to $\zeta_j$ moves with a velocity proportional to $\operatorname{Re}\zeta_j$, is suggested by the Lax-pair time dependence and by soliton resolution \cite{JLPS18cmp}.
\end{remark}

Under soliton resolution \cite{JLPS18cmp}, a solution decomposes for large times into solitons, attached to the eigenvalues, and dispersive radiation. If all eigenvalues lie in the sector $\arg\zeta\le\pi/2-\eta$, every soliton moves with a speed bounded below by a multiple of its frequency, and radiation disperses. This motivates the following spectral strengthening of Conjecture~\ref{conj:A}.

\begin{conjecture}[Spectral form]\label{conj:spectral}
For every $\eta\in(0,\pi/2)$ and $m>0$ there is $C(\eta,m)$ with the following property. Every $q_0\in\Sch$ with $\Mass(q_0)\le m$, whose transmission coefficient has no zeros in the sector $\{\pi/2-\eta<\arg\zeta<\pi/2+\eta\}$, satisfies $\LS_1(\Phi(\cdot)q_0)\le C(\eta,m)\Mass(q_0)$. Moreover, one can take $C(\eta,m)\lesssim_m\eta^{-1}$.
\end{conjecture}

By Remark~\ref{rem:spectral}, Schwartz data of mass $m<2\pi$ satisfy the hypothesis with $\eta=(2\pi-m)/4$. Therefore Conjecture~\ref{conj:spectral}, combined with the density argument in the proof of Corollary~\ref{cor:small}, implies Conjecture~\ref{conj:A}, including the rate $(2\pi-m)^{-1}$, which Theorem~\ref{thm:B} shows to be optimal. The spectral formulation also explains why breathers or multi-soliton bound states should not produce stationary structures below $2\pi$. Indeed, a stationary bound state would need all its eigenvalues on the imaginary axis, and each of them contributes $2\pi$ to \eqref{trace}. This last assertion relies on the heuristic velocity--eigenvalue relation of Remark~\ref{rem:spectral} for multi-soliton states.

\subsection{The sharp dependence on the frequency scale}

The lower bound \eqref{thmA-lower} and the soliton heuristics suggest that the factor $\kappa^2$ in \eqref{thmA} is an artefact of the method.

\begin{conjecture}[Sharp quantitative local smoothing]\label{conj:B}
There is a nondecreasing function $C:(0,\infty)\to(0,\infty)$ such that for all $q_0\in\Sch$ and $T>0$,
\[
\LS_T\big(\Phi(\cdot)q_0\big)\le C\big(\Mass(q_0)\big)\big(1+\kappa(q_0)\,T\big)\Mass(q_0).
\]
\end{conjecture}

The stationary soliton shows that the exponent one cannot be improved. A proof would require local smoothing on windows of size much larger than $\kappa^{-1}$ during times of order one. This seems to need a substitute for the parabolic rescaling that preserves the smallness used in \cite{HKNV}.

\subsection{Global-in-time smoothing for small mass}

The estimate \eqref{small2} holds on parabolic cylinders of every size, which suggests asking whether smoothing for small mass holds globally in time. Some care is needed in the formulation. For the free flow with low-frequency data, $|u(t,0)|^2$ decays like $t^{-1}$. Hence $\int_\R\|Wu(t)\|_{\Hdot}^2dt$ can diverge logarithmically, because $\|Wu(t)\|_{\Hdot}$ is comparable to $\|Wu(t)\|_{L^2}$ for such slowly varying functions, which reflects the commutator between $W$ and $|D|^{1/2}$. The natural global quantity is the Kato--Kenig--Ponce--Vega one.

\begin{question}\label{q:global}
Is there $m_2>0$ such that every $q_0\in L^2$ with $\Mass(q_0)\le m_2$ satisfies
\[
\sup_{h\in\R}\int_\R\int_\R W_h(x)^2\,\big|D_x^{1/2}q(t,x)\big|^2\,dx\,dt\le C\,\Mass(q_0)\ ?
\]
\end{question}

For the free Schr\"odinger group the answer is yes, by the identity of \cite{KPV91}. Small-mass solitons, which by Remark~\ref{rem:fast} do not scatter, appear to satisfy this bound, since their frequency and speed are both of order $\lambda\cot(2\vartheta)$. A positive answer would be a nonlinear Kato smoothing theorem at the critical regularity, in a setting without scattering. A natural tool is the family of microscopic conservation laws of \cite{TangXu21,HKNV} with the spectral parameter coupled to time, in the spirit of an interaction Morawetz estimate. We have not been able to close such an argument. The finite-scale analysis of \cite{BaiMajdoub} may help to understand the linear part of the problem.

\subsection{Strichartz-type quantities}

By \cite[Proposition~1.11]{HKNV}, the rescaled stationary solitons also satisfy $\int_{-1}^1\int_\R|q_n|^6\,dx\,dt\to\infty$. (This is also immediate from scaling, since $\int_\R|u_\lambda(t)|^6dx=\lambda^2\int_\R|q_s(0)|^6dx$ for the rescaled stationary soliton $u_\lambda$.) Our methods yield the localized $L^4_{t,x}$ bound \eqref{small3} for small mass, but not an $L^6_{t,x}$ bound. The latter would require control of $\|W_hq\|_{\Hdot}$ in $L^4_t$ rather than in $L^2_t$, since $\|f\|_{L^6}^6\lesssim\|f\|_{L^2}^2\|f\|_{\Hdot}^4$. We do not know whether $\int_{-1}^1\int_\R|W_hq|^6\,dx\,dt\le C\Mass^3$ holds even for small mass. Its threshold, if it exists, is presumably also $2\pi$.

\subsection{Possible routes to Conjecture~\ref{conj:A} and their obstacles}

A natural approach is by contradiction and concentration compactness, in the spirit of Kenig and Merle \cite{KenigMerle06}. Suppose that $\Cm(m)=\infty$ for some $m<2\pi$, and take a minimal such $m$. One would extract a sequence $q_n$ with $\LS_1(q_n)/\Mass(q_n)\to\infty$. One would then rescale it so that the bad behaviour occurs at unit frequency on time intervals of length $T_n\to\infty$, and try to pass to a limit which is almost periodic modulo symmetries and has mass at most $m$. Three obstacles stand in the way. First, a profile decomposition adapted to the DNLS flow in $L^2$ is not available. Bahouri and Perelman \cite{BahouriPerelman22} use profile decompositions at the $H^{1/2}$ and $H^1$ levels, and the absence of uniform continuity of the flow on bounded sets of $H^s$, $s<1/2$ \cite{BiagioniLinares01}, prevents the usual nonlinear profile approximation. Second, minimality of $m$ and the Lipschitz stability theory that usually comes with it have no counterpart here. Third, the rigidity step would have to exclude non-dispersing almost periodic solutions of mass below $2\pi$. Remark~\ref{rem:spectral} suggests doing so through the spectral data, which are conserved and somehow well behaved for $L^2$ data in the sense of \cite[Corollary~2.1]{BahouriPerelmanT}.

A second approach is to strengthen the bootstrap of \cite[Section~4]{HKNV} directly. There, the quartic and higher terms of the current are controlled by two copies of the local smoothing norm times a small factor $\|q\|_{E^\eps}$. Below $2\pi$, one would need to replace this smallness by coercivity. The Hayashi--Ozawa coercivity \cite{HayashiOzawa92,HayashiOzawa94} is an $H^1$-level statement for the gauged energy, and it does not obviously localize to the frequency-dependent quantities $A(\kappa;q)$. Finding a coercive $L^2$-level substitute would be a significant step, for instance a lower bound on $\operatorname{Re}A(\kappa;q)$ for $\Mass(q)<2\pi$ that is uniform across frequency scales. 

\section{Concluding remarks and further directions}\label{sec:concl}

We have isolated a quantitative question left open by the $L^2$ well-posedness theory of \cite{HKNV}: on which balls of $L^2(\R)$ does local smoothing for DNLS hold uniformly? We proved that it holds for small mass, with a scale-invariant Morrey-type form, and that it fails from mass $2\pi$ on. We also proved that the uniform constant must blow up at least like $(2\pi-m)^{-1}$, which is exactly the rate exhibited by the soliton family. Finally, we made explicit the dependence on the data in the local smoothing estimate of \cite{HKNV}. It is carried by a single frequency scale $\kappa(q_0)$, and it is polynomial with an exponent between one and two.

The results proved here rest on \cite{HKNV} as a black box together with elementary arguments, and their value lies mainly in the questions they make precise. Conjecture~\ref{conj:A} asserts that $2\pi$, the mass of the stationary soliton and the Hayashi--Ozawa coercivity threshold, is sharp. Conjecture~\ref{conj:spectral} gives a spectral form, and Conjecture~\ref{conj:B} the sharp dependence on the frequency scale. Question~\ref{q:global} asks for global-in-time smoothing at small mass. A positive answer would be a nonlinear Kato smoothing theorem at critical regularity, in a setting where small solitons prevent scattering.

Several further directions are natural. The commuting-flow method of \cite{HKNV} applies, at least formally, to the whole Kaup--Newell hierarchy, and the threshold question can be asked for each flow. For odd flows the relevant soliton velocities change, and the threshold may differ from $2\pi$. The models (1.33) and (1.36) of \cite{HKNV} are obtained from \eqref{DNLS} by a Galilean-type modulation and by a gauge transformation, respectively, and they are well-posed in $L^2$ by \cite[Corollary~1.13]{HKNV}. The threshold can be studied for them as well, keeping in mind that gauge transformations do not commute with spatial localization. For the continuum Calogero--Moser model, where $2\pi$ is already a well-posedness threshold \cite{KLV25,KimKimKwon}, it would be interesting to know whether uniform local smoothing holds up to $2\pi$. Finally, the periodic problem in $L^2(\mathbb{T})$ remains open. Since local smoothing is unavailable on the circle, a different mechanism will be needed there.

\bigskip

\section*{Declarations}
\noindent\textbf{Funding.} This research received no funding.

\noindent\textbf{Conflict of interest.} The author declares no conflict of interest.

\noindent\textbf{Data availability.} Data sharing is not applicable to this article, as no datasets were generated or analysed.


\end{document}